\documentclass[a4paper]{article}
\usepackage{amssymb}
\usepackage{amsthm}
\usepackage{amsmath}
\usepackage[cp1250]{inputenc}
\usepackage[boxed]{algorithm2e}
\usepackage{graphicx}

\newtheorem{theorem}{Theorem}[section]
\newtheorem{lemma}[theorem]{Lemma}
\newtheorem{rem}[theorem]{Remark}
\newtheorem{definition}[theorem]{Definition}
\theoremstyle{definition}

\DeclareMathOperator{\IM}{Im}

\DeclareMathOperator{\BALL}{B}
\DeclareMathOperator{\SQ}{Sq}

\DeclareMathOperator{\RE}{Re}
\DeclareMathOperator{\spect}{spect}
\DeclareMathOperator{\DOM}{dom}
\DeclareMathOperator{\MID}{mid}
\DeclareMathOperator{\REST}{r}
\providecommand{\B}[1]{\BALL\left(#1\right)}
\providecommand{\sq}[1]{\SQ\left(#1\right)}

\providecommand{\im}[1]{\IM\left(#1\right)}
\providecommand{\re}[1]{\RE\left(#1\right)}

\providecommand{\norm}[1]{\lVert#1\rVert}
\providecommand{\jednadruga}{\frac{1}{2}}
\providecommand{\dom}[1]{\DOM\left(#1\right)}
\providecommand{\midI}[1]{\MID\left(#1\right)}
\providecommand{\restI}[1]{\REST\left(#1\right)}

\providecommand{\fartail}{T_\mathcal{F}}
\providecommand{\neartail}{T_\mathcal{N}}

\providecommand{\fun}[1]{\FuncSty{#1}}

\providecommand{\decay}{order of polynomial decay}

\def\qed{{\hfill{\vrule height5pt width3pt depth0pt}\medskip}}

\begin{document}
\title{Existence of globally attracting fixed points of viscous Burgers equation with constant forcing. A computer assisted proof}
\author{Jacek~Cyranka\footnote{Research has been supported by National Science Centre grant DEC-2011/01/N/ST6/00995.}\\{\small Institute of Computer Science, Jagiellonian University}\\{\small prof. Stanis\l awa \L ojasiewicza 6, 30-348 Krak\'ow, 
Poland}\\{\small jacek.cyranka@ii.uj.edu.pl}}
\date{{\small \today}}
\maketitle
\begin{abstract}
	We present a computer assisted method for proving the existence of globally attracting fixed points of dissipative PDEs.
	An application to the viscous Burgers equation with periodic boundary conditions and a forcing function constant in time is presented as a case
	study. We establish the existence of a locally attracting fixed point by using rigorous numerics techniques. To prove that the fixed point 
	is, in fact, globally attracting we introduce a technique relying on a construction of an 
	absorbing set, capturing any sufficiently regular initial condition after a finite time. Then the absorbing set is rigorously integrated
	forward in time to verify that any sufficiently regular initial condition is in the basin of attraction of the fixed point.
\end{abstract}
\section{Introduction}
	The field of computer assisted proofs for ordinary differential equations (ODEs) is a quite well established and analysed topic. 
	Still, it seems to us that the development of methods for investigating the dynamics of PDEs by performing rigorous computer assisted 
	proofs is at a pioneering stage. 	
	
	In the present paper we develop a computer aided method which is interesting for two main reasons. First, it provides not only a local, but 
	also a global perspective on the dynamics. Second, it allows to establish results which have not been achieved using known analytical techniques. As a 
	case study we present the forced viscous Burgers equation, where the forcing is constant in time and periodic in space. More specifically, we 
	consider the initial value problem with periodic boundary conditions for the equation
	\begin{equation}
		\label{eq:vBEq}
		u_t+u\cdot u_x-\nu u_{xx}=f(x).
	\end{equation}  	
	In the present paper we deal with the case of non-zero forcing, which is not reducible to a linear PDE by the Hopf-Cole transform anymore. 
	
	To our knowledge, there exist two rigorous numerics methods for studying the non-stationary PDE problem using the Fourier basis. 
	The method of self-consistent bounds, 
	presented in the series of papers \cite{ZM}, \cite{Z2}, \cite{Z3}, \cite{ZAKS}, and the method presented in \cite{AK}. Both of them have been 
	applied to the Kuramoto-Sivashinsky equation. In \cite{KKN} authors obtained some rigorous numerics prototype results for a
	non-stationary PDE problem using the Finite Element basis.
	Related work regarding a rigorous numerics study of the global dynamics of PDE includes \cite{DHMO}, and \cite{MPMW}.
	In \cite{FTKS} the viscous Burgers equation with zero forcing was used as an illustration of a computer 
	aided technique of proving existence of stationary solutions.
	
	It has been shown that \eqref{eq:vBEq} belongs to the class of dissipative PDEs (dPDEs) possessing inertial manifolds \cite{V}. 	
	Using our technique we demonstrate that the global attractor exhibited by \eqref{eq:vBEq} is in fact a unique stable fixed point. 
	In \cite{JKM} it was shown that for any viscosity and the time independent forcing the attractor of \eqref{eq:vBEq} is a single point.  
	This is a stronger than ours result, but the methods in present paper have also some advantages. Contrary to the approach from \cite{JKM} 
	we are not invoking any unconstructive functional analysis techniques, thus the speed of convergence could be obtained from our construction. 
	Moreover, we are not using the maximum principle, so our method should apply to a class of systems of PDEs. 
	
	To establish the existence of an 
	attracting fixed point locally, we use the computer techniques from \cite{ZAKS}. We construct a small neighbourhood of a candidate for the fixed 
	point and prove the existence and uniqueness of a fixed point within said neighbourhood by calculating an explicit upper bound for the logarithmic 
	norm. In case of the negative logarithmic norm, we claim that there exists a locally attracting fixed point. On the other hand, we show the  
	global existence of solutions by constructing trapping regions inspired by  the analogical sets constructed for the Navier-Stokes equations 
	\cite{MS}, \cite{ES}, see also \cite{ZNS}.
 
	We link those results by constructing an absorbing set, which captures any initial condition after a 
	finite time. Then we integrate the obtained absorbing set forward in time rigorously until it is mapped into a small region with the established existence of an attracting fixed point within. By doing so, we verify that any initial condition is in the basin of attraction of the fixed point. 
	The aforementioned elements applied together give an original technique 
	that allows to extend the property of attractiveness obtained locally on a small region to a global fact. We would like to stress that our 
	method concerns the evolution of dPDEs in time, not only the stationary problem. Moreover, it is worth pointing out that we do not restrict ourselves
	by assuming zero spatial average, i.e. $\int_Q{u(t,x)}\,dx=0$ on a domain $Q$, which was often assumed in related work, see for 
	instance \cite{V} or \cite{FTKS}. Our theory can be applied when zero is replaced by any number. We remark that exclusively in the case of non-zero 
	spatial average the equation \eqref{eq:vBEq} admits travelling wave-like solutions. 

	An example result obtained with the presented method is the following
\providecommand{\paperExampleL}{-0.162445}
\providecommand{\paperExampleTotalTime}{3.135}
\providecommand{\paperExampleNrOfSteps}{627}
\providecommand{\paperExampleExecutionTime}{4.36}	
\providecommand{\paperExampleNu}{[2,2.1]}
\providecommand{\paperExampleForcing}{1.6\cos{2x}-2\sin{3x}}
\providecommand{\paperExampleBall}{\sum_{k=1}^{3}{\beta_k\sin{kx}+\gamma_k\cos{kx}},\ \beta_k,\gamma_k\in\left[-0.03,0.03\right]}
\providecommand{\paperExampleAzero}{\pi}
	\begin{theorem}
		\label{thm:main}
		For any $\nu\in\paperExampleNu$ and\\$f\in \left\{x\mapsto\paperExampleForcing+\paperExampleBall\right\}$ there exists a steady state solution of 
		\eqref{eq:vBEq}, which is unique and attracts globally any initial data $u_0$ satisfying $u_0\in C^4$ and 
		$\int_0^{2\pi}{u_0(x)\,dx}=\paperExampleAzero$.
	\end{theorem}
	Other examples are given in Section~\ref{sec:exampleTheorem}. The function $\paperExampleForcing$, added to the forcing was chosen as an example to show 
	that our method is not limited to the simpler case of low energy forcings. Note that Theorem~\ref{thm:main} covers a whole set of forcing 
	functions within a ``ball'' $\paperExampleBall$. To achieve this we used the interval arithmetic in a way to be explained later. 

	By using the presented algorithm we could prove a more general case, namely replace in Theorem~\ref{thm:main} $\beta_k,\ \gamma_k$ with 
	arbitrary continuous functions $\beta_k(t),\ \gamma_k(t)$, such that $\beta_k(t),\ \gamma_k(t)\in[-0.03,0.03]$ for $t\geq 0$. This
	will be exploited in the next paper \cite{CZ} where we prove existence of globally attracting periodic orbits for viscous Burgers equation
	with nonautonomous forcing.

	This paper is dependent on \cite{Z3} and \cite{ZAKS}, we recall only crucial definitions and results from the previous works and 
	focus on the new elements. Proper references are always provided whenever necessary. We are convinced that the presented techniques are applicable 
	to higher dimensional dPDEs, including the Navier-Stokes equations, and we will address this problem in our forthcoming papers.		
	
	We organize the paper as follows: the first part comprises the theory and it is concluded by the proof of Theorem~\ref{thm:main} in 
	Section~\ref{sec:exampleTheorem}. A presentation and discussions of the algorithms follows.
\section{The viscous Burgers equation}
	\label{sec:burgers}
    As the viscous Burgers equation we consider the following PDE\\
    \begin{equation*}
        \frac{\partial u}{\partial t}+u\cdot\frac{\partial u}{\partial x}-\nu\bigtriangleup u=0\quad{\mbox{in}\ } \Omega,\quad t>0,
    \end{equation*}
	where $\nu$ is a positive \textit{viscosity constant}. The equation was proposed by Burgers (1948) as a mathematical model of turbulence. 
	Later on it was successfully showed that the Burgers equation models certain gas dynamics (Lighthill (1956)) and acoustic (Blackstock (1966)) phenomena,
	see e.g. \cite{Wh}. 
	We consider the equation on the real line $\Omega:=\mathbb{R}$ with periodic boundary conditions and \emph{a constant in time forcing $f$}, i.e.
	\begin{equation*}
		u\colon\mathbb{R}\times[0,\mathcal{T})\to\mathbb{R},
	\end{equation*}
	\begin{equation*}
		f\colon\mathbb{R}\to\mathbb{R},
	\end{equation*}
	\begin{subequations}
		\label{eq:burgers}
		\begin{align}
			&u_t+u\cdot u_x-\nu u_{xx}=f(x),\quad x\in\mathbb{R},\ t\in[0, \mathcal{T}),\label{eq:burgers1}\\			
			&u(x,t)=u(x+2k\pi, t),\quad x\in\mathbb{R},\ t\in[0, \mathcal{T}),\ k\in\mathbb{Z},\label{eq:burgers3}\\
			&f(x)=f(x+2k\pi),\quad x\in\mathbb{R},\ k\in\mathbb{Z},\label{eq:burgers4}\\
			&u(x,0)=u_0(x),\quad x\in\mathbb{R}.\label{eq:burgers2}
		\end{align}
	\end{subequations}
\subsection{The viscous Burgers equation in the Fourier basis}
	In this section we rewrite \eqref{eq:burgers1} using \textit{the Fourier basis} of $2\pi$ periodic functions $\{e^{ikx}\}_{k\in\mathbb{Z}}$. From 
	now on we assume that all functions we use are sufficiently regular to be expanded in the Fourier basis and all necessary Fourier series converge. 
	\begin{definition}
		\label{def:FourierModes}
		Let $u\colon\mathbb{R}\to\mathbb{R}$ be a $2\pi$ periodic function.
		We call $\{a_k\}_{k\in\mathbb{Z}}$ {\rm the Fourier modes} of $u$, where $a_k\in\mathbb{C}$ satisfies 
		\begin{equation}
			\label{eq:FourierModes1}
			a_k=\frac{1}{2\pi}\int_0^{2\pi}{u(x)e^{-ikx}}\,dx,
		\end{equation}moreover, the following equality holds
		\begin{equation}
			\label{eq:FourierModes2}
			u(x)=\sum_{k\in\mathbb{Z}}{a_ke^{ikx}},\quad x\in\mathbb{R}.
		\end{equation}
	\end{definition}
	\begin{definition}
	  Let $\nmid\cdot\nmid\colon\mathbb{R}\to\mathbb{R}$ be given by
	  \begin{equation*}
	    \nmid a\nmid:=\left\{\begin{array}{ll}|a|&\text{ if }a\neq 0,\\1&\text{ if }a=0.\end{array}\right.
	  \end{equation*}
	\end{definition}
	\begin{lemma}
		\label{lem:fourierC}
		Let $\gamma>1$. Assume that $|a_k|\leq\frac{M}{\nmid k\nmid^\gamma}$ for $k\in\mathbb{Z}$. If $n\in\mathbb{N}$ is such that $\gamma-n>1$, then the function
		$u(x)=\sum_{k\in\mathbb{Z}}{a_ke^{ikx}}$ belongs to $C^n$. The series
		\begin{equation*}
			\frac{\partial^s u}{\partial x^s}(x)=\sum_{k\in\mathbb{Z}}{a_k\frac{\partial^s}{\partial x^s}e^{ikx}}
		\end{equation*}
		converges uniformly for $0\leq s\leq n$.
	\end{lemma}
	\begin{lemma}
		Let $u_0$ be an initial value for the problem \eqref{eq:burgers} and $f$ be a forcing. Then \eqref{eq:burgers} rewritten in the 
		Fourier basis becomes 
		\begin{subequations}
			\label{eq:burgers_infinite}
			\begin{align}	
				&\frac{d a_k}{d t}=-i\frac{k}{2}\sum_{k_1\in\mathbb{Z}}{a_{k_1}\cdot a_{k-k_1}}+\lambda_k a_k+f_k,\quad k\in\mathbb{Z},\label{eq:burgers_infinite1}\\
				&a_k(0)=\frac{1}{2\pi}\int_0^{2\pi}{u_0(x)e^{-ikx}}\,dx,\quad k\in\mathbb{Z},\label{eq:burgers_infinite2}\\
				&f_k=\frac{1}{2\pi}\int_0^{2\pi}{f(x)e^{-ikx}}\,dx,\quad k\in\mathbb{Z},\label{eq:burgers_infinite3}\\
				&\lambda_k=-\nu k^2.\label{eq:burgers_infinite4}
			\end{align}
		\end{subequations}		
	\end{lemma}
	For the proof refer \cite{SuppMat}.	
	\begin{definition}
		\label{def:symmetricGalerkinProjection}
		For any given number $m>0$ {\rm the $m$-th Galerkin projection} of \eqref{eq:burgers_infinite1} is
		\begin{equation}
			\label{eq:symmetricGalerkinProjection}
			\frac{d a_k}{d t}=-i\frac{k}{2}\sum_{\substack{|k-k_1|\leq m\\|k_1|\leq m}}{a_{k_1}\cdot a_{k-k_1}}+\lambda_k a_k+f_k,\quad|k|\leq m.
		\end{equation}
	\end{definition}	
	Note that in our case $\{a_k\}_{k\in\mathbb{Z}}$ are not independent. The solution $u$ of \eqref{eq:burgers} is real valued, which implies 
	that 
	\begin{equation}
		\label{eq:first}
		a_k=\overline{a_{-k}}.
	\end{equation} 
	Note that condition \eqref{eq:first} is invariant under all Galerkin projections \eqref{eq:symmetricGalerkinProjection} as long as $f_k=\overline{f_{-k}}$.
	
	In Section~\ref{sec:analytic} and Section~\ref{sec:global} we will assume that the initial condition for \eqref{eq:burgers_infinite} 
	satisfies
	\begin{equation}
	  \label{eq:fixedInt}
	  \frac{1}{2\pi}\int_{0}^{2\pi}{u_0(x)\,dx}=\alpha,\quad\text{for a fixed }\alpha\in\mathbb{R}.
	\end{equation}
	We will require additionally that $f_0=0$, and then \eqref{eq:fixedInt} implies that $a_0(t)$ is constant in time, namely
	\begin{equation}
	  \label{eq:fixedA0}
	  a_0=\alpha.
	\end{equation}
	Note that condition \eqref{eq:fixedA0} is invariant under all Galerkin projections \eqref{eq:symmetricGalerkinProjection} as long as $f_0=0$.
\section{Analytic arguments}
\label{sec:analytic}
	In this section we provide some analytic arguments that we use in proving the global existence and regularity results for solutions of 
	\eqref{eq:burgers}.
\providecommand{\gp}{Galerkin projection of \eqref{eq:burgers_infinite1}}
\providecommand{\gps}{Galerkin projections of \eqref{eq:burgers_infinite1}}
\subsection{Energy as Lyapunov function}
    \begin{definition}
		\label{def:energy}
		\emph{Energy} of (\ref{eq:burgers_infinite1}) is given by the formula
    \begin{equation}
			\label{eq:energy}
			E(\{a_k\})=\sum_{k\in\mathbb{Z}}{|a_k|^2}.
		\end{equation}
		Energy of (\ref{eq:burgers_infinite1}) with $a_0$ excluded is given by the formula
		\begin{equation}
		  \label{eq:energyWithoutZero}
		  \mathcal{E}(\{a_k\})=\sum_{k\in\mathbb{Z}\setminus\{0\}}{|a_k|^2}.
		\end{equation}
    \end{definition}
	The following lemma provides an argument for the statement that \textit{the energy} of \eqref{eq:burgers_infinite1} is being absorbed by a ball 
	whose radius depends on the forcing and the viscosity constant. Basing on this argument, later on, we will construct a trapping region for any \gp\ .
	In particular, any trapping region constructed encloses the absorbing ball.
	\begin{lemma}
        \label{lem:energy}
        For any solution of \eqref{eq:burgers_infinite1} or a \gp\ such that $a_{-k}=\overline{a_k}$ the following equality holds
        \begin{equation}
          \label{eq:energyEquality}
            \frac{d E(\{a_k\})}{d t}= -2\nu\sum_{k\in\mathbb{Z}}{k^2|a_k|^2}+\sum_{k\in\mathbb{Z}}{f_{-k}\cdot a_k}+\sum_{k\in\mathbb{Z}}{f_{k}\cdot a_{-k}}.
        \end{equation}
    \end{lemma}
    \paragraph{	\textit{Proof}} Using the symmetry of the index in \eqref{eq:energy} we rewrite 
    \begin{multline*}
        \frac{d E}{d t}=\sum_{k\in\mathbb{Z}}{(\frac{d a_k}{d t}\cdot a_{-k})}+\sum_{k\in\mathbb{Z}}{(\frac{d a_{-k}}{d t}\cdot a_k)}=\sum_{k\in\mathbb{Z}}{-i\frac{k}{2}\sum_{k_1\in\mathbb{Z}}{a_{k_1}\cdot a_{k-k_1}\cdot a_{-k}}}\\
        +\sum_{k\in\mathbb{Z}}{i\frac{k}{2}\sum_{k_1\in\mathbb{Z}}{a_{k_1}\cdot a_{-k-k_1}\cdot a_{k}}}-2\nu\sum_{k\in\mathbb{Z}}{k^2a_k\cdot a_{-k}}+\sum_{k\in\mathbb{Z}}{f_{-k}\cdot a_k}+\sum_{k\in\mathbb{Z}}{f_{k}\cdot a_{-k}}\\
        =\sum_{k\in\mathbb{Z}}{-ik\sum_{k_1\in\mathbb{Z}}{a_{k_1}\cdot a_{k-k_1}\cdot a_{-k}}}-2\nu\sum_{k\in\mathbb{Z}}{k^2a_k\cdot a_{-k}}+\sum_{k\in\mathbb{Z}}{f_{-k}\cdot a_k}+\sum_{k\in\mathbb{Z}}{f_{k}\cdot a_{-k}}.
    \end{multline*}
    We want to show that $\sum_{|k|\leq N}{k\sum_{\substack{|k_1|\leq N\\|k-k_1|\leq N}}{a_{k_1}\cdot a_{k-k_1}\cdot a_{-k}}}=0$.
    In order to facilitate the proof explanation we denote $S_{N,k}:=\sum_{\substack{|l|\leq N\\|k-l|\leq N}}{a_{k-l}\cdot a_l}$ and 
    $S_N:=\sum_{|k|\leq N}{k\sum_{\substack{|k_1|\leq N\\|k-k_1|\leq N}}{a_{k_1}\cdot a_{k-k_1}\cdot a_{-k}}}=\sum_{|k|\leq N}{kS_{N,k}a_{-k}}$.
    \paragraph{}We proceed by induction, firstly we check if for $N=1$ the thesis is fulfilled
	\begin{equation*}
        S_1=-1(a_{-1}\cdot a_0\cdot a_1+a_0\cdot a_{-1}\cdot a_1)+1(a_1\cdot a_0\cdot a_{-1}+a_0\cdot a_1\cdot a_{-1})=0.		
    \end{equation*}
    We verify the induction step $S_{N-1}=0 \Rightarrow S_N=0$
    \begin{equation*}		
        S_N=S_{N-1}+\sum\left\{\begin{array}{lll}
        a_N\cdot a_{-N+k}\cdot a_{-k}2k&, 0<k<N,&\ \left(S_{I}\right)\\
        a_{-N}\cdot a_{k+N}\cdot a_{-k}2k&, -N<k<0,&\ \left(S_{II}\right)\\
		    S_{N,-N}\cdot a_N(-N)&, k=-N,&\ \left(S_{III}\right)\\
		    S_{N,N}\cdot a_{-N}N&, k=N,&\ \left(S_{IV}\right)        
		\end{array}\right.
    \end{equation*}
    we match elements with the same modes from $\left(S_I\right)$ and $\left(S_{III}\right)$. Let $e(N)=1$ for $N$ even and $e(N)=0$ for $N$ odd, 
    \begin{multline*}\sum_{0<k<N}{a_N\cdot a_{k-N}\cdot a_{-k}(2k-N)}\\
		=\sum_{0<k<\frac{N}{2}}{a_N\cdot a_{k-N}\cdot a_{-k}(2k-N+N-2k)}+e(N)a_N\cdot a^2_{-\frac{N}{2}}(N-N)=0.
    \end{multline*}
    When elements with the same modes from $\left(S_{II}\right)$ and $\left(S_{IV}\right)$ are matched analogously as above the result 
    is also zero.    
    After substitution all that is left is
    \begin{equation*}
        S_N=S_{N-1}+2Na_N\cdot a_0\cdot a_{-N}-2Na_{-N}\cdot a_0\cdot a_N=0.\quad\qed
    \end{equation*}
\subsection{A trapping region for \eqref{eq:burgers_infinite1}}
    \label{sec:analyticalTR}
    
	In this section we provide a forward invariant set for each \gp, called \textit{the trapping region} . 
	If we consider an arbitrary initial condition that is inside a trapping region, then the corresponding trajectory remains in this set in 
	the future. This is an argument for the existence of solutions of each \gp\ 
	within a trapping region. Moreover, due to the existence of a trapping region, the solution of \eqref{eq:burgers_infinite1}, obtained by passing to 
	the limit, conserves the initial regularity. We use this fact to argue that a solution of \eqref{eq:burgers_infinite} with sufficiently regular 
	initial data exists for all times, is unique, and is a classical solution of \eqref{eq:burgers}. Calculations performed in this section were 
	inspired by the trapping regions built for the Navier-Stokes equations, see \cite{MS} and \cite{ES}. 
	\paragraph{\textit{Notation}} Let $l^2(\mathbb{Z})=\left\{\left\{a_k\right\}_{k\in\mathbb{Z}}\colon\sum{|a_k|^2<\infty}\right\}$, where $a_k\in\mathbb{C}$
	for $k\in\mathbb{Z}$. In the sequel the space $l^2(\mathbb{Z})$ will be denoted by $H$. We equip $H$ with the standard scalar 
	product. Let $m>0$, we define $P_m(H)$ to be $\mathbb{C}^{2m+1}$.
	\paragraph{}Formally an element of a Galerkin projection \eqref{eq:symmetricGalerkinProjection} is a finite sequence. In the sequel we will use the
	following embedding, and with some abuse of notation we will use the same symbol to denote the element of infinite dimensional space $H$
	\begin{equation*}
	  P_m(H)\ni\left\{a_{-m},\dots,a_0,\dots,a_m\right\}\equiv\left\{\dots,0,\dots,0,a_{-m},\dots,a_0,\dots,a_m,0,\dots,0,\dots\right\}\in H.
	\end{equation*}
	In consequence we assume the inclusion $P_m(W)\subset W$, for all $W\subset H$.
	\begin{lemma}
		\label{lem:estimateNk}
		Let $\{a_k\}_{k\in\mathbb{Z}}\in H$, $N_k:=-i\frac{k}{2}\sum_{k_1\in\mathbb{Z}}{a_{k_1}\cdot a_{k-k_1}}$. 
		Assume that there exists $C>0$ and $s>0.5$ such that 
		$\{a_k\}_{k\in\mathbb{Z}}$ satisfy $|a_k|\leq\frac{C}{\nmid k\nmid^s}$, $k\in\mathbb{Z}$.\\ Then 
		\begin{equation*}
			|N_k|\leq\frac{\sqrt{E(\{a_k\})}C\left(2^{s-\jednadruga}+\frac{2^{s-1}}{\sqrt{2s-1}}\right)}{|k|^{s-\frac{3}{2}}},\quad k\in\mathbb{Z}\setminus\{0\}.
	    	\end{equation*}
	\end{lemma}
	\paragraph{\textit{Proof}}	
  In order to prove the bound for $N_k$, we split $N_k=N_k^I+N_k^{II}$, and bound $N_k^I$ and $N_k^{II}$ separately 
	\paragraph{Case 1} First, we bound the following sum $N_k^{I}=-i\frac{k}{2}\sum_{k_1}{a_{k_1}\cdot a_{k-k_1}}$,\quad where $|k_1|\leq\frac{1}{2}|k|$\\
    \begin{multline*}
        |N_k^{I}|\leq\sum_{|k_1|\leq\frac{1}{2}|k|}{\jednadruga|k||a_{k_1}||a_{k-k_1}|}\leq\sum_{|k_1|\leq\jednadruga|k|}{\jednadruga|k|\frac{C}{|k-k_1|^{s}}|a_{k_1}|}\\
        \leq\frac{2^{s-1}C}{|k|^{s-1}}\sqrt{\sum_{|k_1|\leq\frac{1}{2}|k|}{|a_{k_1}|^2}}\sqrt{\sum_{|k_1|\leq\frac{1}{2}|k|}{1}}\leq\frac{2^{s-1}\sqrt{2}\sqrt{E(\{a_k\})}C}{|k|^{s-\frac{3}{2}}}.
    \end{multline*}
    \paragraph{Case 2} Second, we bound the remaining part $N_k^{II}=-i\frac{k}{2}\sum_{k_1}{a_{k_1}\cdot a_{k-k_1}}$,\quad where $|k_1|>\jednadruga|k|$\\
    \begin{multline*}
        |N_k^{II}|\leq\sum_{|k_1|>\jednadruga|k|}{\frac{1}{2}|k||a_{k_1}||a_{k-k_1}|}\leq\frac{1}{2}|k|C\sum_{|k_1|>\jednadruga|k|}{\frac{1}{|k_1|^s}|a_{k-k_1}|}\\
		\leq\frac{1}{2}|k|C\sqrt{\sum_{|k_1|>\jednadruga|k|}\frac{1}{|k_1|^{2s}}}\sqrt{\sum_{|k_1|>\jednadruga|k|}{|a_{k-k_1}|^2}}\\
		\leq\frac{1}{2}|k|\sqrt{E(\{a_k\})}C\sqrt{\frac{2^{2s}}{(2s-1)|k|^{2s-1}}}=\frac{\sqrt{E(\{a_k\})}C\frac{2^{s-1}}{\sqrt{2s-1}}}{|k|^{s-\frac{3}{2}}}.
    \end{multline*}
    We used the following estimation due to the convexity
    \begin{equation*}
        \sum_{|k_1|>\jednadruga|k|}{\frac{1}{|k_1|^{2s}}}<2\int_{\jednadruga|k|}^{\infty}\frac{1}{r^{2s}}\,d r=2\left[-\frac{1}{(2s-1)r^{2s-1}}\right]_{\jednadruga|k|}^{\infty}=\frac{2^{2s}}{(2s-1)|k|^{2s-1}}.
    \end{equation*}    
    After summing together Case~1 and Case~2
    \begin{equation*}		
        |N_k|\leq|N_k^{I}|+|N_k^{II}|=\frac{\sqrt{E(\{a_k\})}C\left(2^{s-\jednadruga}+\frac{2^{s-1}}{\sqrt{2s-1}}\right)}{|k|^{s-\frac{3}{2}}}.
    \end{equation*}
    holds for any $k\in\mathbb{Z}\setminus\{0\}$.\qed
\providecommand{\D}{2^{s-\jednadruga}+\frac{2^{s-1}}{\sqrt{2s-1}}}
\providecommand{\Ef}{E\left(\left\{f_k\right\}\right)}
\providecommand{\Ezero}{\mathcal{E}}
\providecommand{\Fassumption}{$f_k=\overline{f_{-k}}$, $f_k=0$ for $|k|>J$ and $f_0=0$}
\providecommand{\satisfiesFassumption}{$\{f_k\}$ satisfies \Fassumption}
\providecommand{\energyBound}{\widetilde{\mathcal{E}}}
\providecommand{\invariantSubspace}{restricted to the invariant subspace given by $a_k=\overline{a_{-k}}$ and $a_0=\alpha$}
	\begin{theorem}
    \label{thm:analyticTrappingRegion}
    Let $\{a_k\}_{k\in\mathbb{Z}}\in H$, $\alpha\in\mathbb{R}$, $J>0$, $s>0.5$, $E_0=\frac{\Ef}{\nu^2}$, $\energyBound>E_0$, $D=\D$, 
    $C>\sqrt{\energyBound}N^s$, $N>\max{\left\{J,\left(\frac{\sqrt{\energyBound+\alpha^2}D}{\nu}\right)^{2}\right\}}$. Assume that
	  \satisfiesFassumption. Then 
	\begin{equation*}
		W_0(\energyBound, N, C, s,\alpha)=\{\{a_k\}\ |\ \Ezero(\{a_k\})\leq \energyBound,\ |a_k|\leq\frac{C}{|k|^s}\text{ for }|k|>N\}
	\end{equation*}
    is a trapping region for each \gp\ \invariantSubspace.
    \end{theorem}
	\paragraph{	\textit{Proof}}
	We first show that 
	\begin{equation}
		\label{eq:energyDecreasing}
		\text{if }\Ezero(\{a_k\})>E_0=\frac{\Ef}{\nu^2}\text{ then }\frac{d \Ezero(\{a_k\})}{dt}<0.
	\end{equation}
	Under the assumption $f_0=0$ we have 
	$\sum_{k\in\mathbb{Z}}{|f_{-k}||a_k|} + \sum_{k\in\mathbb{Z}}{|f_k||a_{-k}|}=\sum_{k\in\mathbb{Z}\setminus\{0\}}{|f_{-k}||a_k|} + \sum_{k\in\mathbb{Z}\setminus\{0\}}{|f_k||a_{-k}|}$, and 
	$\frac{d a_0}{d t}=0$, the latter implies that 
  $\frac{d E}{d t}=\frac{d\Ezero}{d t}$.
  
  Taking the square root of $\Ezero(\{a_k\})>\frac{\Ef}{\nu^2}$ gives
  \begin{equation*}
    \nu\sqrt{\Ezero(\{a_k\})}>\sqrt{\Ef},
  \end{equation*}
  multiplying both of the sides by $2\sqrt{\sum_{k\in\mathbb{Z}\setminus\{0\}}{|a_k|^2}}$ gives
  \begin{equation*}
    2\nu\sum_{k\in\mathbb{Z}\setminus\{0\}}{|a_k|^2}>2\sqrt{\sum_{k\in\mathbb{Z}\setminus\{0\}}{|f_k|^2}}\sqrt{\sum_{k\in\mathbb{Z}\setminus\{0\}}{|a_k|^2}},    
  \end{equation*}
  moreover the following inequalities are satisfied
	\begin{equation*}
		2\nu\sum_{k\in\mathbb{Z}}{k^2a_k\cdot a_{-k}}\geq 2\nu\sum_{k\in\mathbb{Z}\setminus\{0\}}{|a_k|^2}>2\sqrt{\sum_{k\in\mathbb{Z}\setminus\{0\}}{|f_k|^2}}\sqrt{\sum_{k\in\mathbb{Z}\setminus\{0\}}{|a_k|^2}}\geq\sum_{k\in\mathbb{Z}\setminus\{0\}}{|f_{-k}||a_k|} + \sum_{k\in\mathbb{Z}\setminus\{0\}}{|f_k||a_{-k}|}.
	\end{equation*}
	Simply, the linear term dominates the forcing term in \eqref{eq:energyEquality}, i.e. 
	\begin{equation}
	  \label{eq:linearDominates}
	  2\nu\sum_{k\in\mathbb{Z}}{k^2a_k\cdot a_{-k}}>\sum_{k\in\mathbb{Z}\setminus\{0\}}{|f_{-k}||a_k|} + \sum_{k\in\mathbb{Z}\setminus\{0\}}{|f_k||a_{-k}|}.
	\end{equation}
	The condition \eqref{eq:linearDominates} is sufficient to satisfy \eqref{eq:energyDecreasing}.
	
	Next observe that the condition $|a_k|\leq\frac{C}{|k|^s}$ is satisfied for all 
	$\{a_k\}\in W_0$ and $k\in\mathbb{Z}\setminus\{0\}$. 
	Since $\Ezero(\{a_k\})\leq \energyBound$ and $|a_k|\leq\sqrt{\energyBound}$ for $k\neq 0$,
    \begin{equation*}
        |a_k|\leq \sqrt{\energyBound}\leq\frac{C}{|k|^s}\text{ because }C>\sqrt{\energyBound}N^s.
    \end{equation*}
	\paragraph{}Now, we shall check if the vector field points inwards on $\partial W_0$.
  For $\{a_k\}\in\partial W_0$ such that $\Ezero(\{a_k\})=\energyBound$ and $\energyBound>\frac{\Ef}{\nu^2}$ vector field points inwards from \eqref{eq:energyDecreasing}. 
	Let us pick a point $\{a_k\}_{k\in\mathbb{Z}}\in\partial W$ such that $|a_k|=\frac{C}{|k|^s}$ for some $|k|>N$, and perform calculations to check if the diminution condition $\frac{d |a_k|}{d t}<0$
	holds. Observe that $E(\{a_k\}_{k\in\mathbb{Z}})\leq\energyBound+\alpha^2$ and we apply Lemma~\ref{lem:estimateNk} with 
	$E(\{a_k\}_{k\in\mathbb{Z}})$ replaced by $\energyBound+\alpha^2$.
    \begin{gather*}
    	\frac{d |a_k|}{d t}<-\nu|k|^2\frac{C}{|k|^s}+\frac{D\sqrt{\energyBound+\alpha^2}C}{|k|^{s-\frac{3}{2}}}<0,\\
        \nu|k|^2\frac{C}{|k|^s}>\frac{D\sqrt{\energyBound+\alpha^2}C}{|k|^{s-\frac{3}{2}}},\\
        \nu\sqrt{|k|}>D\sqrt{\energyBound+\alpha^2},\\
		|k|>\left(\frac{D\sqrt{\energyBound+\alpha^2}}{\nu}\right)^2,
    \end{gather*}
    $\frac{d |a_k|}{d t}<0$ holds if $|k|>\frac{D^2\left(\energyBound+\alpha^2\right)}{\nu^2}$. The proof is complete because 
    $|k|>N>\frac{D^2\left(\energyBound+\alpha^2\right)}{\nu^2}$.
	$\qed$
\section{Global results}
	\label{sec:global}
\providecommand{\subspaceH}{\overline{H}}
	\begin{definition}
	  The subspace $\subspaceH\subset H$ is defined by
	  \begin{equation*}
      \subspaceH:=\left\{\{a_k\}\in H\colon\text{there exists }0\leq C<\infty\text{ such that }|a_k|\leq\frac{C}{\nmid k\nmid^4}\text{ for }k\in\mathbb{Z}\right\}.	  
	  \end{equation*}
	\end{definition}
	\paragraph{\textit{Notation}} Let $l>0$, we define $P_l(H)$ to be $\mathbb{C}^{2l+1}$. From now on by $\varphi^l(t, x)$ we denote the solution of $l$-th \gp\,
	at a time $t>0$, with an initial value $x\in P_l(H)$. By $\{a_k\}_{|k|\leq l}$ we denote an initial condition $x\in P_l(H)$. The operator 
	$N_k$ is the nonlinear part of \eqref{eq:burgers_infinite1}, and is defined by 
	$N_k(\{a_k\}_{k\in\mathbb{Z}}):=-i\frac{k}{2}\sum_{k_1\in\mathbb{Z}}{a_{k_1}\cdot a_{k-k_1}}$ for $k\in\mathbb{Z}$. For a sequence of complex
	numbers $\{c_k\}_{k\in\mathbb{Z}}$ let $c_{k,j}$ denotes the $j$-th component of $c_k$ for $k\in\mathbb{Z}$ and $j=1,2$, complex numbers 
	are considered as elements of $\mathbb{R}^2$ here.
	
	Let $P_l(H)\ni\{a^l_k(t)\}_{|k|\leq l}:=\varphi^l(t, \{a_k\}_{|k|\leq l})$, $t>0$, $l>0$. Observe that $\{a^l_k(t)\}_{|k|\leq l}$ is well defined, 
	as solutions for each Galerkin projection of \eqref{eq:burgers_infinite1} exist for all times $t>0$ due to Theorem~\ref{thm:analyticTrappingRegion} 
	(existence of a trapping region) and are unique due to the fact that \eqref{eq:symmetricGalerkinProjection} is a finite system of ODEs with a 
	locally Lipschitz right-hand side. We will drop the index $l$ when it is known either from the context or irrelevant in the context.				
	
\providecommand{\Wassumption}{$H\supset W$ be a trapping region for $l$-th \gp\ \invariantSubspace\ for all $l>M_1$}
\providecommand{\WassumptionSimple}{$H\supset W$ be a trapping region for $l$-th \gp\ for all $l>M_1$}
\providecommand{\aassumption}{$a_k=\overline{a_{-k}}$ and $a_0=\alpha$}
	\begin{lemma}
		\label{lem:energyDissipation}
		Let $\alpha\in\mathbb{R}$, $J>0$, $M_1\geq 0$, $E_0=\frac{\Ef}{\nu^2}$, $\energyBound>E_0$. Assume that \satisfiesFassumption. 
		Let \Wassumption. \\
		There exists a finite time $t_1=t_1(W)\geq 0$ such that $\Ezero\left(\varphi^l(t_1,P_l(\{a_k\}_{k\in\mathbb{Z}}))\right)\leq \energyBound$ 
		holds uniformly for all $\{a_k\}_{k\in\mathbb{Z}}\in W$ and $l>M_1$.
	\end{lemma}
	\paragraph{	\textit{Proof}} Let us take $\{\hat{a}_k\}_{k\in\mathbb{Z}}$ from the boundary of $W$, such that $\Ezero(\{a_k\}_{k\in\mathbb{Z}})\leq\Ezero(\{\hat{a}_k\}_{k\in\mathbb{Z}})$
	for all $\{a_k\}_{k\in\mathbb{Z}}\in W$.
	Let $\Ezero\left(\left\{\hat{a}_k\right\}_{k\in\mathbb{Z}}\right)=\Ezero_I$ be 
	the initial energy. It is enough to take either $t_1=0$ if $\Ezero_I\leq\energyBound$ or 
	$t_1(W)=\frac{1}{2\nu\varepsilon}\ln{\frac{\Ezero_I}{\energyBound}}$ if $\Ezero_I>\energyBound$, where $\varepsilon=\left(1-\sqrt{\frac{E_0}{\energyBound}}\right)$.\\
	To see this, we calculate in a similar fashion as in the proof of Theorem~\ref{thm:analyticTrappingRegion}. Let $\Ezero_I>\energyBound$, by 
	Lemma~\ref{lem:energy} and the assumption that $f_0=0$ (observe that in this case $\frac{d E}{d t}=\frac{d\Ezero}{d t}$, because $a_0$ is 
	a constant) we have
	\begin{equation*}
		\frac{d\Ezero}{dt}\leq -2\nu \Ezero+2\sqrt{\Ezero}\sqrt{E(\{f\})}=-2\nu \Ezero\left(1-\frac{\sqrt{E(\{f\})}}{\nu\sqrt{\Ezero}}\right)\leq
		-2\nu \Ezero\left(1-\sqrt{\frac{E_0}{\energyBound}}\right),
	\end{equation*}
	therefore by Gronwall's inequality
	\begin{equation*}
		\Ezero(t)\leq e^{-2\nu t\left(1-\sqrt{\frac{E_0}{\energyBound}}\right)}\Ezero_I.
	\end{equation*}
	We set $t_1=t$, where $t$ satisfies $e^{-2\nu\varepsilon t}\Ezero_I=\energyBound$. The time 
	$t_1$ is uniform for the trapping region $W$, because $\Ezero_I$ is the maximal energy within the trapping region $W$.\qed
\providecommand{\numbersInRR}{The numbers $f_k\in\mathbb{C}$, $N_k\in\mathbb{C}$ and $a_k\in\mathbb{C}$ are considered as elements of $\mathbb{R}^2$}
	\begin{lemma}
		\label{lem:akbk1}
		Let $M_1\geq 0$, $k\in\mathbb{Z}\setminus\{0\}$, $j=1,2$, $\lambda_k$ be the $k$-th eigenvalue \eqref{eq:burgers_infinite4}, \WassumptionSimple.
		\numbersInRR.\\
		Assume that $N_{k,j}^\pm\in\mathbb{R}^2$ are bounds such that
		\begin{equation*}
			\left(N_{k,1}(\{a_k\}_{k\in\mathbb{Z}}), N_{k,2}(\{a_k\}_{k\in\mathbb{Z}})\right)\in[N_{k,1}^-, N_{k,1}^+]\times[N_{k,2}^-, N_{k,2}^+]\text{ for all }\{a_k\}_{k\in\mathbb{Z}}\in W.
		\end{equation*}
		Then for any $\varepsilon>0$ there exists a finite time $\hat{t}>0$ such that for all $l>\max{\{M_1, |k|\}}$ and $t\geq\hat{t}$
		$a_k^l(t)$ with any initial condition in $P_l(W)$ satisfies
		\begin{equation*}
			\left(a_{k,1}^l(t), a_{k,2}^l(t)\right)\in \left[b_{k,1}^-,b_{k,1}^+\right]\times\left[b_{k,2}^-,b_{k,2}^+\right]+\left[-\varepsilon,\varepsilon\right]^2,
		\end{equation*}
		where $b_{k,j}^\pm=\frac{N_{k,j}^\pm+f_{k,j}}{-\lambda_k}$.
	\end{lemma}
	\paragraph{	\textit{Proof}} In the calculations we drop the index $l$ denoting the Galerkin projection dimension and the index $j$ denoting the coordinate 
	for better clarification, for instance instead of $a_{k,j}^l$ we write $a_k$. We perform the calculations for the 
	first and the second component simultaneously; thus, we finally obtain two values $t_{k, 1}>0$ and $t_{k,2}>0$. For any \gp\ from 
	$\frac{d a_k}{dt}\leq\lambda_k a_k+N_k^++f_k^+$, $\frac{d a_k}{dt}\geq\lambda_k a_k+N_k^-+f_k^-$ it follows that 
	\begin{equation*}
		a_k(t)\geq\left(a_k^--b_k^-\right)e^{\lambda_k t}+b_k^-,\quad a_k(t)\leq\left(a_k^+-b_k^+\right)e^{\lambda_k t}+b_k^+,	
	\end{equation*}
	where $a_k^\pm$ are bounds such that $\left(a_{k,1}, a_{k,2}\right)\in[a_{k,1}^-, a_{k,1}^+]\times[a_{k,2}^-, a_{k,2}^+]$, which exist as 
	the initial condition is contained in a compact trapping region. \\
	Because $\lambda_kt=-\nu k^2t<0$ for any $t>0$ ($k\in\mathbb{Z}\setminus\{0\}$ by assumption) it follows that for a sufficiently large time $t_k>0$ 
	we have $\left(\left|a_k^+-b_k^+\right|+\left|a_k^--b_k^-\right|\right)e^{\lambda_k t} \leq \varepsilon$ for any $t\geq t_k$. It is enough to take 
	\begin{equation*}
		t_k=-\ln{\frac{\varepsilon}{\left(\left|a_k^+-b_k^+\right|+\left|a_k^--b_k^-\right|\right)}}/{\nu k^2}.
	\end{equation*}
	Finally $\hat{t}:=\max\left\{t_{k,1}, t_{k,2}\right\}$.\qed
	\begin{lemma}
		\label{lem:akbk2}
		Let $J>0$, $M_1\geq 0$, \WassumptionSimple. 
		Assume that $C_a$, $s_a$ are numbers such that 
		\begin{equation*}
		  |a_k|\leq\frac{C_a}{|k|^{s_a}}\text{ for }|k|>M_1\text{, and for all }\{a_k\}\in W.
		\end{equation*}
		Assume that $\{f_k\}$ satisfies $f_k=0$ for $|k|>J$, $f_0=0$, and 
		$C_N$, $s_N$ are numbers such that
		\begin{equation*}
			\left|N_k(\{a_k\}_{k\in\mathbb{Z}})\right|\leq\frac{C_N}{|k|^{s_N}}\text{ for }|k|>M_1.
		\end{equation*}
		\\Then for any $\varepsilon>0$ there exists a finite time $\hat{t}\geq 0$ such that for all $l>M_1$ and $t\geq\hat{t}$
		 $\left\{a_k^l(t)\right\}_{|k|\leq l}$ with any initial condition in $P_l(W)$ satisfy 
		\begin{equation*}
				|a_k^l(t)|\leq\frac{C_b+\varepsilon}{|k|^{s_b}}\text{ for }|k|>M_1,
		\end{equation*}
		where $C_b=\left(C_N+\max_{0<|k|\leq J}{\left\{|f_k||k|^{s_N}\right\}}\right)/\nu$, $s_b=s_N+2$.
	\end{lemma}
	\paragraph{	\textit{Proof}} 
	We will use the same notation as in Lemma~\ref{lem:akbk1}. For any \gp\  from the fact that 
	$\frac{d a_k^l}{d t}\leq\lambda_k\left(a_k^l+\frac{N_k^++f_k^+}{\lambda_k}\right)$ and 
	$\frac{d a_k^l}{d t}\geq\lambda_k\left(a_k^l+\frac{N_k^-+f_k^-}{\lambda_k}\right)$, it follows
	\begin{gather*}
	  a_k^l(t)\leq\left(\frac{C_a}{|k|^{s_a}}-\frac{C_b}{|k|^{s_b}}\right)e^{\lambda_{k}t}+\frac{C_b}{|k|^{s_b}},\\
	  a_k^l(t)\geq\left(-\frac{C_a}{|k|^{s_a}}+\frac{C_b}{|k|^{s_b}}\right)e^{\lambda_{k}t}-\frac{C_b}{|k|^{s_b}},
	\end{gather*}
	for $|k|>M_1$. Due to the fact that $s_b>s_a$
	\begin{equation*}
		|a_k^l(t)|\leq\frac{C_a(k_{max}(t))^{s_b-s_a}e^{\lambda_{k_{max}(t)}t}+C_b}{|k|^{s_b}},\quad |k|>M_1
	\end{equation*}
	for all $l>M_1$ and $t>0$, where $C_b=\left(C_N+\max_{0<|k|\leq J}{\left\{|f_k||k|^{s_N}\right\}}\right)/\nu$, $s_b=s_N+2$, $k_{max}(t)$ is the value 
	at which the maximum of $f_t(k)=C_ak^{s_b-s_a}e^{\lambda_k t}$ is attained. Analogically, for a 
	sufficiently large time $t_\mathcal{F}>0$ 
	\begin{equation*}
		{C_a(k_{max}(t))^{s_b-s_a}e^{\lambda_{k_{max}(t)}t}\leq\varepsilon},\quad t\geq t_\mathcal{F},
	\end{equation*} therefore
	\begin{equation*}
		|a_k^l(t)|\leq\frac{C_b+\varepsilon}{|k|^{s_b}},\quad t\geq t_\mathcal{F},\ l>M_1.
	\end{equation*}
	Finally, the obtained time $t_\mathcal{F}$ is uniform with respect to the projection dimension $l$.\qed
	\begin{lemma}
		\label{lem:energyEstimate}
		Let $\hat{E}>0$. The following estimate holds
		\begin{equation*}
			\left|N_k(\{a_k\}_{k\in\mathbb{Z}})\right|\leq\jednadruga|k|\hat{E}
		\end{equation*}
		for all $\{a_k\}_{k\in\mathbb{Z}}\in\left\{\{a_k\}_{k\in\mathbb{Z}}\in H\ |\ E(\{a_k\}_{k\in\mathbb{Z}})\leq \hat{E}\right\}$.
	\end{lemma}
	\paragraph{	\textit{Proof}} Let $\{a_k\}_{k\in\mathbb{Z}}\in \left\{\{a_k\}_{k\in\mathbb{Z}}\in H\ |\ E(\{a_k\}_{k\in\mathbb{Z}})\leq \hat{E}\right\}$.	
	We start with the easy estimate $\left|N_k(\{a_k\}_{k\in\mathbb{Z}})\right|\leq\jednadruga|k|\sum_{k_1\in\mathbb{Z}}{|a_{k-k_1}||a_{k_1}|}$, 
	by the Cauchy-Schwarz inequality 
	$\left|N_k(\{a_k\}_{k\in\mathbb{Z}})\right|\leq{\jednadruga|k|\sqrt{\sum_{k_1\in\mathbb{Z}}{|a_{k_1}|^2}}\sqrt{\sum_{k_1\in\mathbb{Z}}{|a_{k-k_1}|^2}}}$,
	which is the following energy estimate $\left|N_k(\{a_k\}_{k\in\mathbb{Z}})\right|\leq\jednadruga|k|\hat{E}$.\qed
	\paragraph{}Now, we shall introduce \textit{the absorbing sets}. For any initial condition there exists a finite time after which the solutions of Galerkin 
	projections are trapped in an absorbing set. We use absorbing sets as a tool for studying the global dynamics of 
	\eqref{eq:burgers_infinite1}.
	\begin{definition}
		\label{def:absorbingSet}
		Let $M_1>0$. A set $A\subset H$ is called {\rm the absorbing set} for large \gps, if for any initial 
		condition $\{a_k\}_{k\in\mathbb{Z}}\in H$ 
		there exists a finite time $t_1\geq 0$ such that for all $l>M_1$ and $t\geq t_1$ 
		$\varphi^l\left(t, P_l(\{a_k\}_{k\in\mathbb{Z}})\right)\in P_l(A)$.
	\end{definition}
	In what follows we will often call the absorbing set for large \gps\ simply \emph{the absorbing set}.
	
	In the next result, to show the existence of an absorbing set, we construct analytically an absorbing set. Furthermore, we construct 
	absorbing sets with any \decay. Later on, 
	in the context of a computer assisted proof of the main theorem, we will construct an absorbing set using \textit{the interval arithmetic}. 
	Accomplishing this task requires the established existence of an absorbing set with a sufficiently large \decay. 
	\begin{lemma}
		\label{lem:absorbingSet}
		Let $\alpha\in\mathbb{R}$, $\varepsilon>0$, $J>0$, $M_1\geq 0$, $E_0=\frac{\Ef}{\nu^2}$, $\energyBound>E_0$. Assume that \satisfiesFassumption. Put
		\begin{eqnarray*}
		  s_i&=&i/2\text{ for }i\geq 2,\\
		  D_i&=&2^{s_i-\frac{1}{2}}+\frac{2^{s_i-1}}{\sqrt{2s_i-1}}\text{ for }i\geq 2,\\
		  C_2&=&\varepsilon+\left(\jednadruga\left(\energyBound+\alpha^2\right)+\max_{0<|k|\leq J}{\frac{|f_k|}{|k|}}\right)/\nu,\\
		  C_i&=&\varepsilon+\left(C_{i-1}\sqrt{\energyBound+\alpha^2}D_{i-1}+\max_{0<|k|\leq J}{|k|^{s_i-2}|f_k|}\right)/\nu\text{ for }i>2.		  		  
		\end{eqnarray*}\\		
		Then for all $i\geq 2$
		\begin{equation*}
			H\supset W_i\bigl(\energyBound,M_1,\varepsilon,\alpha\bigr):=\left\{\{a_k\}_{k\in\mathbb{Z}}\ |\ \Ezero(\{a_k\}_{k\in\mathbb{Z}})\leq\energyBound,\ |a_k|\leq\frac{C_i}{|k|^{s_i}}\text{ for }|k|>M_1\right\}
		\end{equation*} 
		is an absorbing set for large \gps\ \invariantSubspace.		 
	\end{lemma}
	\paragraph{	\textit{Proof}} Let $\energyBound>E_0$, $\{\widehat{a_k}\}_{k\in\mathbb{Z}}$ be an arbitrary initial condition for 
	\eqref{eq:burgers_infinite}, $E^{max}:=\max\left\{\Ezero\left(\{\widehat{a_k}\}_{k\in\mathbb{Z}}\right), \energyBound\right\}$. Let 
	$C_0\geq 0$ and $s_0> 0$ be constants such that 
	\begin{equation}
		\label{eq:W_0}
		W_0:=\left\{\{a_k\}_{k\in\mathbb{Z}}\ |\ \Ezero(\{a_k\}_{k\in\mathbb{Z}})\leq E^{max},\ |a_k|\leq\frac{C_0}{\nmid k\nmid^{s_0}}\right\}
	\end{equation} 
	is a trapping region for each Galerkin projection of \eqref{eq:burgers_infinite1} enclosing $\{\widehat{a_k}\}_{k\in\mathbb{Z}}$. This trapping 
	region exists due to Theorem~\ref{thm:analyticTrappingRegion}. Note that a trapping region can be scaled to make it enclose an arbitrary sufficiently
	smooth initial condition. It follows from Lemma~\ref{lem:energyDissipation} that there exists a finite time $t_1\geq 0$ such that for all $\{a_k\}_{k\in\mathbb{Z}}\in W_0$ 
	and $l>M_1$
	\begin{equation}
		\label{eq:energyW0}
		\Ezero(\varphi^l(t_1,P_l(\{a_k\}_{k\in\mathbb{Z}})))\leq \energyBound.
	\end{equation}
	We define $W_1:=W_0\cap\left\{\{a_k\}_{k\in\mathbb{Z}}\ |\ \Ezero(\{a_k\}_{k\in\mathbb{Z}})\leq \energyBound\right\}$. 
	From \eqref{eq:energyW0} and that $W_0$, $W_1$ are trapping regions we immediately have that $\varphi^l(t, P_l(\{a_k\}_{k\in\mathbb{Z}})\in W_1$
	for all $\{a_k\}_{k\in\mathbb{Z}}\in W_0$, $t\geq t_1$ and $l>M_1$. Using Lemma~\ref{lem:energyEstimate} we bound the nonlinear part  
	\begin{equation}
		\label{eq:energyEstimate}
		|N_k(\{a_k\}_{k\in\mathbb{Z}})|\leq\jednadruga|k|\left(\energyBound+\alpha^2\right)\text{ for all }\{a_k\}_{k\in\mathbb{Z}}\in W_1.
	\end{equation} 
	It follows from Lemma~\ref{lem:akbk2} that there exists a finite time $t_2\geq t_1$ such that for all 
	$t\geq t_2$ and $l>M_1$, $\{a_k^l(t)\}_{|k|\leq l}$ with any initial condition in $P_l(W_1)$ satisfy
	\begin{equation}
		\label{eq:C_2}
		|a_k^l(t)|\leq\frac{C_2}{|k|}\text{ for }|k|>M_1.
	\end{equation}
  It is important to start with 
	the energy estimate \eqref{eq:energyEstimate} to bound the nonlinear part $N_k$ because the goal is to estimate $|a_k|$ uniformly with 
	respect to $C_0$ and $s_0$ \eqref{eq:W_0}. We emphasize that $C_2$ from \eqref{eq:C_2} does not depend on $C_0$ and $s_0$. Having the bound 
	\eqref{eq:C_2}, we construct the following absorbing set 
	\begin{equation*}
		W_2:=\left\{\{a_k\}_{k\in\mathbb{Z}}\ |\ \Ezero(\{a_k\}_{k\in\mathbb{Z}})\leq\energyBound,\ |a_k|\leq\frac{C_2}{|k|}\text{, for }|k|>M_1\right\}.
	\end{equation*}
	Due to Lemma~\ref{lem:estimateNk} the following estimate holds
	\begin{equation*}
		\left|N_k(\{a_k\}_{k\in\mathbb{Z}})\right|\leq\frac{C_2\sqrt{\energyBound+\alpha^2}D_2}{|k|^{-\jednadruga}},
	\end{equation*}
	for all $\{a_k\}_{k\in\mathbb{Z}}\in W_2$. Due to Lemma~\ref{lem:akbk2} again there exists a finite time $t_3\geq t_2$ such that for all $t\geq t_3$ 
	and $l>M_1$, $\{a_k^l(t)\}_{|k|\leq l}$ with any initial condition in $P_l(W_2)$ satisfy
	\begin{equation*}
		|a_k^l(t)|\leq\frac{C_3}{|k|^{\frac{3}{2}}}\text{ for }|k|>M_1,
	\end{equation*}
	where $C_3=\varepsilon+\left(C_2\sqrt{\energyBound+\alpha^2}D_2+\max_{0<|k|\leq J}{\frac{|f_k|}{|k|^{\jednadruga}}}\right)/\nu$. Having this bound, we construct the following 
	absorbing set
	\begin{equation*}
		W_3:=\left\{\{a_k\}_{k\in\mathbb{Z}}\ |\ \Ezero(\{a_k\}_{k\in\mathbb{Z}})\leq \energyBound,\ |a_k|\leq\frac{C_3}{|k|^{\frac{3}{2}}}\text{ for }|k|>M_1\right\}.
	\end{equation*}
	Note the gain of $\jednadruga$ in the \decay\ of $\{a_k\}_{k\in\mathbb{Z}}$ in $W_3$ compared to $W_2$. From applying Lemma~\ref{lem:estimateNk} 
	and Lemma~\ref{lem:akbk2} further we obtain a sequence of times $t_3<t_4<\dots<t_n<\dots$ such that
	\begin{equation*}
		|a_k^l(t_3)|\leq\frac{C_3}{|k|^{s_3}},\dots,\ |a_k^l(t_n)|\leq\frac{C_n}{|k|^{s_n}},\dots,\text{ for }|k|>M_1,
	\end{equation*}
	with $s_i=i/2$, $C_i=\varepsilon+\left(C_{i-1}\sqrt{\energyBound+\alpha^2}D_{i-1}+\max_{0<|k|\leq J}{|k|^{s_i-2}|f_k|}\right)/\nu$.
	The obtained $W_i$, $i\geq 2$ are absorbing sets for large \gps, which follows from the construction and that $C_i$ for all $i\geq 2$ depend on the 
	energy $\energyBound$ and $\alpha$ only.\qed
	\begin{rem}
		Assume the same as in Lemma~\ref{lem:absorbingSet}. The inclusion $W_i\subset\subspaceH$ holds for all $i\geq 8$, where $W_i$ is
		an absorbing set proved to exist in Lemma~\ref{lem:absorbingSet}.
	\end{rem}
	\begin{lemma}
		\label{lem:refinementOfAbsorbingSet}
		Let $k\in\mathbb{Z}\setminus\{0\}$, $\varepsilon>0$, $\lambda_k$ denotes the $k$-th eigenvalue \eqref{eq:burgers_infinite4}. 
		Let $H\supset A$ be an absorbing set for large \gps. 
		\numbersInRR. Assume that $N_k^\pm\in\mathbb{R}^2$ are bounds such that		
		\begin{equation*}
		  \left(N_{k,1}(\{a_k\}_{k\in\mathbb{Z}}), N_{k,2}(\{a_k\}_{k\in\mathbb{Z}})\right)\in[N_{k,1}^-, N_{k,1}^+]\times[N_{k,2}^-, N_{k,2}^+]\text{ for all }\{a_k\}_{k\in\mathbb{Z}}\in A.
		\end{equation*}
		Then $A\cap\left\{\{a_k\}_{k\in\mathbb{Z}}\ |\ a_k\in\left[b_{k,1}^-,b_{k,1}^+\right]\times\left[b_{k,2}^-,b_{k,2}^+\right]+\left[-\varepsilon,\varepsilon\right]^2\right\}$ 
		is also an absorbing set for large \gps, where $b_{k,j}^\pm=\frac{N_{k,j}^\pm+f_{k,j}}{-\lambda_k}$, $j=1,2$.
	\end{lemma}	
	\paragraph{	\textit{Proof}} Immediate consequence of Lemma~\ref{lem:akbk1}.\qed	
\section{General method of self-consistent bounds.}
\label{sec:generalMethod}
	The same symbols as in the preceding part are used in a more general context. For the purpose of the presented work we call a dissipative PDE a PDE of the following type
		\begin{equation}	
			\label{eq:dPDE}
			\frac{du}{dt}=Lu+N(u,Du,\dots,D^r u)+f=F(u),
		\end{equation}
	where $u(x,t)\in\mathbb{R}^n$, $x\in\mathbb{T}^d$, ($\mathbb{T}^d=(\mathbb{R}/2\pi)^d$ is a $d$-dimensional torus), $L$ is a 
	linear operator, $N$ a polynomial and by $D^su$ we denote the collection of $s$-th order partial derivatives of $u$. The right-hand side 
	contains a constant in time forcing function $f$. We require that $L$ is diagonal in the Fourier basis $\{e^{ikx}\}_{k\in\mathbb{Z}^d}$
	\begin{equation*}
		Le^{ikx}=\lambda_k e^{ikx}
	\end{equation*}
	and the \textit{eigenvalues} $\lambda_k$ satisfy
	\begin{subequations}
		\label{eq:dPDEassumptions}
		\begin{gather}
			\lambda_k=-\nu(|k|)|k|^p,\\
			0<\nu_0\leq\nu(|k|)\leq \nu_1,\quad\text{for }|k|>K_-,\\
			p>r,
		\end{gather}
		for some $v_0>0$, $v_1\geq v_0$ and $K_-\geq 0$, $r$ is the maximal order of derivatives appearing in the nonlinear part \eqref{eq:dPDE},
		$|\cdot|$ is the Euclidean norm.
	\end{subequations}	
\subsection{Self-consistent bounds}
	We recall, in the context of dPDEs, the definition of self-consistent bounds from \cite{Z3}. 
	Let $H$ be a Hilbert space, actually $L^2$ or one of its 
	subspaces in the context of dPDEs. We assume that a domain of $F$, the right hand side of \eqref{eq:dPDE}, is dense in $H$. By a solution of \eqref{eq:dPDE} 
	we understand a function $u\colon[0,\mathcal{T})\to\dom{F}$ such that $u$ is differentiable and \eqref{eq:dPDE} is satisfied 
	for all $t\in[0,\mathcal{T})$ and $\mathcal{T}$ is a maximal time of the existence of solution. We assume that there is a set $I\subset\mathbb{Z}^d$ and a 
	sequence of subspaces $H_k\subset H$ for $k\in I$ such that $\dim{H_k}=d_1 <\infty$, $H_k$ and $H_{k'}$ are mutually orthogonal for $k\neq k'$
	and $H=\overline{\oplus_{k\in I}{H_k}}$. Let $A_k\colon H\to H_k$ be the orthogonal projection onto $H_k$, for each $u\in H$ holds 
	$u=\sum_{k\in I}{u_k}=\sum_{k\in I}{A_k u}$. Analogously if $B$ is a function with the range in $H$, then $B_k(u)=A_kB(u).$ 
	
	We assume that a a metric space $(T,\rho)$ is provided, for $X\subset T$ by $\overline{X}$ we 
	denote the closure of $X$, by $\partial X$ we denote the boundary of $X$. For $n>0$ we set 
	$X_n=\oplus_{|k|\leq n,k\in I}{H_k}$, $Y_n=X_n^\perp$. By $P_n\colon H\to X_n$ and $Q_n\colon H\to Y_n$ we denote the orthogonal 
	projections onto $X_n$ and $Y_n$ respectively, $T\supset B(c,r)=\left\{x\in T\colon \rho(c,x)<r\right\}$ denotes a ball with the centre at
	$c$ and the radius $r$.

	\begin{definition}{\rm \cite[Def. 2.1]{Z3}}
		We say that $F\colon H\supset\dom{F}\to H$ is admissible if the following conditions are satisfied for any $n>0$ such that 
		$\dim{X_n}>0$
		\begin{itemize}
			\item $X_n\subset\dom{F}$
			\item $P_nF\colon X_n\to X_n$ is a $C^1$ function
		\end{itemize}
	\end{definition}
	\begin{definition}{\rm \cite[Def. 2.3]{Z3}} 
		\label{def:scs1}
		Assume $F$ is an admissible function. 
		Let $m,\ M\,\in\mathbb{R}$ with $m\leq M$. Consider an object consisting of: a compact set $W\subset X_m$ 
		and a sequence of compact sets $B_k\subset H_k$ for $|k|>m,\ k\in I$. We define the conditions \textbf{C1}, \textbf{C2}, \textbf{C3}, 
		\textbf{C4a} as follows:
		\paragraph{\textbf{C1}} For $|k|>M,\ k\in I$ holds $0\in B_k$.
		\paragraph{\textbf{C2}} Let $\hat{a}_k\colon=\max_{a\in B_k}{\norm{a}}$ for $|k|>m,\ k\in I$ and then $\sum_{|k|>m,k\in I}{\hat{a}^2_k}<\infty$.\\
		In particular
		\begin{equation*}
			W\oplus\Pi_{|k|>m}{B_k}\subset H
		\end{equation*}
		and for every $u\in W\oplus\Pi_{k\in I,|k|>m}{B_k}$ holds , $\norm{Q_n u}^2\leq\sum_{|k|>n,k\in I}{\hat{a}_k^2}$.
		\paragraph{\textbf{C3}} The function $u\mapsto F(u)$ is continuous on $W\oplus\Pi_{k\in I,|k|>m}{B_k}\subset H$.\\
		Moreover, if we define for $k\in I$, $\hat{f}_k=\max_{u\in W\oplus\Pi_{k\in I,|k|>m}{B_k}}{\left|F_k(u)\right|}$, then
		$\sum{\hat{f}_k^2}<\infty$.
		\paragraph{\textbf{C4a}} For $|k|>m,\ k\in I\ B_k$ is given by \eqref{def:ball1} or \eqref{def:ball2}
		\begin{align}
			B_k &= \overline{B(c_k,r_k)},\quad r_k>0,\label{def:ball1}\\
			B_k &= \Pi_{s=1}^{d_1}\left[a_{k,s}^-,a_{k,s}^+\right],\quad a_{k,s}^-<a_{k,s}^+\label{def:ball2}.
		\end{align}
		Let $u\in W\oplus\Pi_{|k|>m}{B_k}$, $F_{k,s}$ be the $s$-th component of $F_k$. Then for $|k|>m$ holds:
		\begin{itemize}
			\item if $B_k$ is given by \eqref{def:ball1} then
			\begin{equation*}
				u_k\in\partial_{H_k}{B_k}\Rightarrow\left(u_k-c_k|F_k(u)\right)<0.
			\end{equation*}
			\item if $B_k$ is given by \eqref{def:ball2} then
			\begin{align*}
				u_{k,s}=a_{k,s}^-&\Rightarrow F_{k,s}(u)>0,\\
				u_{k,s}=a_{k,s}^+&\Rightarrow F_{k,s}(u)<0.
			\end{align*}
		\end{itemize}
	\end{definition}
	\begin{definition}{\rm\cite[Def. 2.4]{Z3}}
		\label{def:scs2}
		Assume $F$ is an admissible function. 
		Let $m,M\in\mathbb{R}$ with $m\leq M$. Consider an object consisting of: a compact set $W\subset X_m$ and
		a sequence of compacts $B_k\subset H_k$ for $|k|>m,k\in I$. We say that set $W\oplus\Pi_{k\in I,|k|>m}{B_k}$ forms {\rm self-consistent bounds for}
		$F$ if conditions \textbf{C1}, \textbf{C2}, \textbf{C3} are satisfied.\\
		If additionally condition \textbf{C4a} holds, then we say that $W\oplus\Pi_{k\in I,|k|>m}{B_k}$ forms {\rm topologically self-consistent bounds for }$F$. 
	\end{definition}				
	\paragraph{}We start our approach by replacing a sufficiently regular $u$ and $f$ in \eqref{eq:dPDE} by the Fourier series, i.e. 
	$u(x,t)=\sum_{k\in\mathbb{Z}^d}{a_k(t)e^{ikx}}$ and $f(x)=\sum_{k\in\mathbb{Z}^d}{f_ke^{ikx}}$. 
	We obtain a system of ODEs describing the evolution of the coefficients $\{a_k\}_{k\in\mathbb{Z}^d}$, 
	where $a_k$ is the coefficient corresponding to $e^{ikx}$
	\begin{equation}
		\label{eq:evolutionFourierModes}
		\frac{d a_k}{d t}=F_k(a)=L_k(a)+N_k(a)+f_k=\lambda_k a_k+N_k(a)+f_k ,\quad k\in\mathbb{Z}^d.
	\end{equation}
	The method works for dPDEs only. The Burgers equation on the real line with forcing, which is the subject of the case study given in this 
	paper is in fact a dPDE.
	\begin{lemma}
	\label{lem:burgersDPDE}
	Let $\nu$ be the viscosity constant in \eqref{eq:burgers1}, then \eqref{eq:burgers1} satisfies the conditions \eqref{eq:dPDEassumptions}
	with $d=1$, $r=1$, $p=2$, $\nu(k)=\nu$, $\lambda_k=-\nu k^2$.
	\end{lemma}

	In our approach we solve the system of equations \eqref{eq:evolutionFourierModes} instead of \eqref{eq:dPDE}. \eqref{eq:evolutionFourierModes}
	is defined on $l^2=\left\{\{a_k\}\colon\,\sum{|a_k|^2<\infty}\right\}$ space or one of its subspaces. We associate $a_k$ with the coefficient
	corresponding to $e^{ikx}$ in the Fourier expansion of $u$.
	Assuming that the initial condition $u_0\in H$ is sufficiently regular, then \eqref{eq:dPDE} and \eqref{eq:evolutionFourierModes} are equivalent. 
	In our approach we expand $u_0$ in the Fourier basis to get the initial value for all the variables $\left\{a_k(0)\right\}_{k\in\mathbb{Z^d}}$. 
	We argue that the solution of \eqref{eq:evolutionFourierModes} 
	is defined for all times $t>0$. Moreover, the solution conserves its initial regularity due to the existence of \textit{trapping regions} and is, in fact, a classical solution of \eqref{eq:dPDE}. For the details refer to Section~\ref{sec:local} and Section~\ref{sec:exampleTheorem}.  
	
	To establish the notation in the next sections we provide	
	\begin{definition}
		\label{def:tail}
		Given an object $W\oplus\Pi_{|k|>m}{B_k}$, $W\subset X_m$ and a sequence of compact sets $B_k\subset H_k$ for $|k|>m$, $m,M\in\mathbb{R}_{+}$,
		$m\leq M$
		\begin{itemize}
			\item $W$ is called {\rm the finite part},
			\item $\Pi_{|k|>m}{B_k}$ is called {\rm the tail} and denoted by $T$,
			\item $\Pi_{m<|k|\leq M}{B_k}$ is called {\rm the near tail} and denoted by $\neartail$,
			\item $\Pi_{|k|> M}{B_k}$ is called {\rm the far tail} and denoted by $\fartail$.
		\end{itemize}
	\end{definition}
	$\neartail$ is the finite part of a tail, whereas $\fartail$ is the infinite part of a tail. In fact in our approach we use $\fartail$ of the form
	\begin{equation}
		\label{eq:farTail}
		\fartail:=\prod_{|k|>M}{\overline{B(0, C/|k|^s)}},\quad C\in\mathbb{R}_+,\ s\geq d+p+1.
	\end{equation}
	\paragraph{}First of all, any $F$ in \eqref{eq:evolutionFourierModes} is admissible, because any finite truncation of a $l^2$ series is in 
	the domain of $F$, and the Galerkin projection of the right-hand side, being a smooth function, is a polynomial.
  $W\oplus T\subset H$ with $\fartail$ defined in \eqref{eq:farTail} 
	satisfies conditions C1, C2 and C3 of Definition~\ref{def:scs1} with $I=\mathbb{Z}^d$, 
	in particular $F$ in \eqref{eq:evolutionFourierModes} is a continuous function on $W\oplus T$. This property was proved in 
	\cite[Theorem 3.6]{Z3}, i.e. $W\oplus T$ forms self-consistent bounds for 
	\eqref{eq:evolutionFourierModes} and equivalently forms a self-consistent bounds for \eqref{eq:dPDE}.
	It is allowable to associate the finite part $W$ with the near tail $\neartail$, but we keep the distinction because of the different 
	treatment of both in the algorithm.
	
	We do not address here the question if solutions of a general dPDE \eqref{eq:dPDE} exist and are unique as it was thoroughly answered 
	in \cite{Z3}, see \cite[Theorem 3.7]{Z3}. 
\section{Local existence and uniqueness}			
	\label{sec:local}
		\paragraph{}Regarding local existence and uniqueness we rely on results from \cite{ZAKS}. For the sake of completeness we recall the main 
		theorems. The same symbols as in the preceding part are used in a more general context.
		\begin{definition}{\rm\cite[Def. 3.1]{ZAKS}}
		  A decomposition of $H$, into into a sum of subspaces is called a block decomposition of $H$ if the following conditions are satisfied.
      \begin{enumerate}
        \item $H=\oplus_i{H_i}$, 
        \item for every $i$ $h_i=dim\,H_i\leq h_{max}<\infty$,
        \item for every $i$ $H_i=\langle e_{i_1},e_{i_2},\dots,e_{i_{h_i}}\rangle$,
        \item if $dim\,H=\infty$, then there exists $k$ such that for $i>k$ $h_i=1$.
      \end{enumerate}
		\end{definition}
		\paragraph{Notation} In this section we adopt the notation from \cite{ZAKS}, 
		namely, we make a distinction between blocks and one dimensional spaces spanned by $\langle e_i\rangle$. For the blocks we use 
		$H_{(i)}=\langle e_{i_1},\dots,e_{i_k}\rangle$, where $(i)=(i_1,\dots,i_k)$. The symbol $H_i$ will always mean the subspace generated by 
		$e_i$. For a block decomposition of $H$ and block $(i)$, we set $\dim{(i)}=\dim{H_{(i)}}$. For any $x\in H$ by $x_{(i)}$ we will denote
		a projection of $x$ onto $H_{(i)}$, by $P_{(i)}$ we will denote an orthogonal projection onto $H_{(i)}$. 
		For $x\in\mathbb{R}^n$ we set $|x|$ to be the Euclidean norm.
		We define the norm (\emph{the 
		block-infinity norm}) by $|x|_{b,\infty}=\max_{(i)}{|P_{(i)}x|}$.
		
		For any norm $||\cdot||$ on $\mathbb{R}^n$ we use the notion of \emph{the logarithmic norm} of a matrix.
		\begin{definition}{\rm\cite[Definition 3.4]{ZAKS}}
		  Let , $Q$ be a square matrix, then we call
		  \begin{equation*}
		    \mu(Q)=\limsup_{h>0,h\to 0}{\frac{||I+hQ||-1}{h}}
		  \end{equation*}
  	  \emph{the logarithmic norm} of $Q$.
		\end{definition}
		\begin{definition}
			Let $R\subset H$, $R$ is convex, $l>0$, $x\in X_l$, $\varphi^l(t, x)$ be the local flow inducted by the $l$-th Galerkin projection of 
			\eqref{eq:evolutionFourierModes}. We call 
			$P_l(R)$ {\rm a trapping region} for the $l$-th Galerkin projection of \eqref{eq:evolutionFourierModes} if 
			$\varphi^l(t, P_l(R))\subset P_l(R)$ for all $t>0$ or equivalently the vector field on the boundary of $P_l(R)$ points inwards.
		\end{definition}		
		\begin{theorem}{\rm\cite[Thm. 3.7]{ZAKS}}
		\label{thm:7}
		Assume that $R\subset H$, $R$ is convex and $F$ satisfies conditions C1, C2, C3. Assume that we have a block decomposition of $H$, such
		that condition \textbf{Db} holds
		\paragraph{Db} there exists $l\in\mathbb{R}$ such that for any (i) and $x\in R$
		\begin{equation}
		\label{eq:logarithmic}
			\mu\left(\frac{\partial F_{(i)}}{\partial x_{(i)}}\left(x\right)\right)+\sum_{(k),(k)\neq(i)}{\left|\frac{\partial F_{(i)}}{\partial x_{(k)}}\left(x\right)\right|}\leq l
		\end{equation}
		Assume that $P_n(R)$ is a trapping region for the $n$-dimensional Galerkin projection of \eqref{eq:evolutionFourierModes} for all $n>M_1$. 
		Then
		\begin{enumerate}
			\item \textbf{Uniform convergence and existence}. For a fixed $x_0\in R$, let $x_n\colon[0,\infty]\to P_n(R)$ be a solution of 
			$x^{\prime}=P_n(F(x)),\ x(0)=P_nx_0$. Then $x_n$ converges uniformly in a max-infinity norm on compact intervals to a function
			$x^{*}\colon[0,\infty]\to R$, which is a solution of \eqref{eq:evolutionFourierModes} and $x^{*}(0)=x_0$. The convergence of $x_n$ on compact 
			time intervals is uniform with respect to $x_0\in R$.
			\item \textbf{Uniqueness within $R$}. There exists only one solution of the initial value problem \eqref{eq:evolutionFourierModes}, $x(0)=x_0$ for 
			any $x_0\in R$ such that $x(t)\in R$ for $t>0$.
			\item \textbf{Lipschitz constant}. Let $x\colon[0,\infty]\to R$ and $y\colon[0,\infty]\to R$ be solutions of \eqref{eq:evolutionFourierModes}, then
			\begin{equation*}
				\left|y(t)-x(t)\right|_{b,\infty}\leq e^{lt}\left|x(0)-y(0)\right|_{b,\infty}
			\end{equation*}
			\item \textbf{Semidynamical system}. The map $\varphi\colon\mathbb{R}_+\times R\to R$, where $\varphi(\cdot,x_0)$ is a unique solution
			of the equation \eqref{eq:evolutionFourierModes} such that $\varphi(0,x_0)=x_0$ defines a semidynamical system on $R$, namely
			\begin{itemize}
				\item $\varphi$ is continuous,
				\item $\varphi(0,x)=x$,
				\item $\varphi(t, \varphi(s, x))=\varphi(t+s,x)$.
			\end{itemize}
		\end{enumerate}	
		\end{theorem}
		\paragraph{}The following Theorem is the main tool used to prove the existence of a locally attracting fixed point.
		\begin{theorem}{\rm\cite[Thm. 3.8]{ZAKS}}
		\label{thm:8}
		The same assumptions on $R, F$ and a block decomposition $H$ as in Theorem~\ref{thm:7}. Assume that $l<0$. 
		\paragraph{} Then there exists a fixed point for \eqref{eq:evolutionFourierModes} $x^*\in R$, unique in $R$, such that for every $y\in R$
		\begin{equation*}
			\begin{split}				
				&\left|\varphi(t,y)-x^*\right|_{b,\infty}\leq e^{lt}\left|y-x^*\right|_{b,\infty},\quad\text{for }t\geq 0,\\
				&\lim_{t\to\infty}{\varphi(t,y)}=x^*.
			\end{split}
		\end{equation*}
		\end{theorem}
\section{Proof of Theorem~\ref{thm:main}}
\label{sec:exampleTheorem}
	Now we are ready to prove Theorem~\ref{thm:main}. The complete algorithm that we used to prove Theorem~\ref{thm:main} and other results 
	\eqref{eq:table} are demonstrated in Section~\ref{sec:cap}. This proof is a prototype for any other result that is obtained using the algorithm, 
	however each case requires construction of different sets. The sets and all the relevant numbers used in the proof of Theorem~\ref{thm:main} are 
	presented in Appendix~\ref{sec:proofData}. 
	\paragraph{\textit{Proof of Theorem~\ref{thm:main}}} Let $u_0\in C^4$ be an arbitrary initial condition satisfying $\int_0^{2\pi}{u_0(x)\,dx}=\paperExampleAzero$, 
	$\left\{a_k\right\}_{k\in\mathbb{Z}}$ be the Fourier coefficients of $u_0$, i.e. $u_0=\sum{a_ke^{ikx}}$. Let $A\subset\subspaceH$ be an absorbing 
	set for large \gps, which exists due to Lemma~\ref{lem:absorbingSet} (for instance $W_8$). Firstly, the existence of a locally attracting fixed point 
	for \eqref{eq:burgers_infinite} is established by constructing a set $\widetilde{W\oplus T}\subset\subspaceH$ 
	satisfying the assumptions of Theorem~\ref{thm:8}, using the interval arithmetic. This is constructed in step~\ref{step:trappingRegion}
	of Algorithm from Section~\ref{sec:cap}. Observe that $\widetilde{W\oplus T}$ is a trapping region for $m$-th \gp\ for all $m>\widehat{m}$ and 
	the logarithmic norm \eqref{eq:logarithmic} is bounded from above by $l<0$. The purpose of the notation $\widetilde{W\oplus T}$ is to keep 
	the consistency with the description of the algorithm used for proving this theorem in Section~\ref{sec:cap}. $\widetilde{W\oplus T}$ satisfies the conditions C1, 
	C2 and C3 of Definition~\ref{def:scs1}, i.e. forms self-consistent bounds for \eqref{eq:burgers_infinite1}, the conditions 
	C1 and C2 are satisfied trivially, because $\overline{H}\subset H$. The condition C3 is also satisfied, the 
	right-hand side of \eqref{eq:burgers_infinite1}, denoted here by $F$, is continuous on $\widetilde{W\oplus T}$. First, notice that 
	$\widetilde{W\oplus T}\subset\dom{F}$ because of the following inequality, let $u\in\widetilde{W\oplus T}$
  \begin{equation}
    \label{eq:Fl2}
    |F(u)_k|\leq\frac{D_1}{|k|^{5/2}}+\frac{\nu D_2}{|k|^2}\leq\frac{\tilde{D}}{|k|^2},\quad |k|>\widehat{m},
  \end{equation} 
  therefore $F(u)\in H$. The continuity of $F$ on $\widetilde{W\oplus T}$ follows from the general theorem \cite[Theorem 3.7]{Z3}. All the 
  assumptions of \cite[Theorem 3.7]{Z3} are satisfied here,
  i.e. \eqref{eq:burgers_infinite1} belongs to the proper class, see Lemma~\ref{lem:burgersDPDE}, and the order of decay of 
  $\widetilde{W\oplus T}$ is sufficient, see \eqref{eq:Fl2}. 
  
  By Theorem~\ref{thm:8} within $\widetilde{W\oplus T}$ there exists a locally 
  attracting fixed point for \eqref{eq:burgers_infinite1}. 		
	
 	Then, $V\oplus\mathnormal{\Theta}\subset\subspaceH$, an absorbing set for large \gps\ satisfying 
 	\begin{equation}
 		\label{eq:rigorousIntegration}
 		\varphi^m\left(t,P_m(V\oplus\mathnormal{\Theta})\right)\subset P_m\left(\widetilde{W\oplus T}\right),
 	\end{equation}
 	for all $t\geq \widehat{t}$ and $m>\widehat{m}$, is constructed. This is constructed in Algorithm from Section \ref{sec:scb}.
 	The absorbing set $V\oplus\mathnormal{\Theta}$ forms self-consistent bounds for \eqref{eq:burgers_infinite1} and thus
 	\eqref{eq:rigorousIntegration} is verified by rigorous integration of $V\oplus\mathnormal{\Theta}$ forward in time using Algorithm~\ref{alg:main}
 	presented hereafter. From \eqref{eq:rigorousIntegration} and the fact that $V\oplus\mathnormal{\Theta}$ is an absorbing set for large \gps\ it follows that
	\begin{equation*}
		\varphi^m\left(t, P_m\left(\left\{a_k\right\}_{k\in\mathbb{Z}}\right)\right)\in P_m\left(\widetilde{W\oplus T}\right).
	\end{equation*} 
	after a finite time and for all $m>\widehat{m}$. Therefore $\left\{a_k\right\}_{k\in\mathbb{Z}}$ is located in the basin of attraction of 
	the fixed point for \eqref{eq:burgers_infinite1}. The sets $\widetilde{W\oplus T}$ and $V\oplus\mathnormal{\Theta}$ are presented in 
	Appendix~\ref{sec:proofData}.
	
	To close the proof we will argue that the fixed point for \eqref{eq:burgers_infinite1} is the steady state solution of \eqref{eq:burgers1}. 
	There exists $C>0$ such that the Fourier coefficients $\left\{a_k\right\}_{k\in\mathbb{Z}}$ of $u_0$ satisfy 	
	\begin{equation}
		\label{eq:regularity}
		|a_k|\leq\frac{C}{\nmid k\nmid^4}.
	\end{equation}
	Let $W_0\subset\subspaceH$ be a trapping region enclosing $\{a_k\}_{k\in\mathbb{Z}}$, and let $\{a_k(t)\}_{k\in\mathbb{Z}}$ be the unique solution of 
	\eqref{eq:burgers_infinite} existing for all times $t>0$, $\{a_k\}_{k\in\mathbb{Z}}\in W_0$ due to Theorem~\ref{thm:7}. The solution is unique, 
	as the logarithmic norm on $W_0$ is bounded, see e.g. \cite{ZNS}. Moreover, the solution conserves the initial regularity \eqref{eq:regularity}. 
	The sequence $\{a_k(t)\}_{k\in\mathbb{Z}}$ for $t>0$ is a 
	classical solution of \eqref{eq:burgers}, as from Lemma~\ref{lem:fourierC}, the condition \eqref{eq:regularity} suffices to $\sum{a_k e^{ikx}}$ 
	and every term that appears in \eqref{eq:burgers1} converge uniformly. Therefore, the solution of \eqref{eq:burgers_infinite} within 
	$\widetilde{W\oplus T}$ is in fact the classical solution of \eqref{eq:burgers}, in particular, the fixed point of \eqref{eq:burgers_infinite1} is 
	the steady state solution of \eqref{eq:burgers1}.\qed
	\paragraph{}In the table below we present example results which we obtained using our algorithm.	
	
	\begin{table}[h!]
	\begin{equation*}
		\label{eq:table}
		\begin{array}{|c|c|c|c|c|c|c|c|c|c|c|}
		\hline
		\mathbf{\nu} & \mathbf{\int_0^{2\pi}{u_0(x)\,dx}}    & \mathbf{E_0}   & \mathbf{\varepsilon}   &   \mathbf{m} &\eqref{eq:logarithmic}\ \mathbf{l< }&     \textbf{1.}           	  &   \textbf{2.}        & \textbf{3.} & \textbf{4.} & \textbf{5.}\\\hline\hline
		   [10,10.1]  &       14\pi       &  	0.5 	  &         0.001		   &	5		  &   -9.94489     &      20.06		   &       1041       	  & \checkmark  & \checkmark  & \checkmark \\\hline		
		    [4,4.1]   &       4\pi        &  	0.5		  &         0.001		   &	7		  &   -2.65147     &      61.41		     &      1305     	 	  & \checkmark  & \checkmark  & \checkmark \\\hline
	\paperExampleNu   &		\paperExampleAzero			  &  	0.82 	  &     	 0.03   	   &	3	  	  &\paperExampleL  &\paperExampleTotalTime	   &\paperExampleNrOfSteps& \checkmark  & \checkmark  & \checkmark \\\hline 
		      1       &       0.4\pi         &  	0.25	  &    	0.0001     &	20		  &  -0.0442416   &     452.23       &      452       	  & \checkmark  & \checkmark  & \checkmark \\\hline
		      0.5     &       0.1\pi 		     &  	0.08	  &    	0.0001		 &	20		  &  -0.0456      &     556.73       &      629       	  & \checkmark  & \checkmark  & \checkmark \\\hline	
		     0.15     &       0              &  	0.22	  &       0 		   &	40  	  &  1340.95      &      26.19                &         -       	  & \checkmark  &   &            \\\hline
 		\end{array}
	\end{equation*}
	\caption{Data from example results}
	\label{table}
	\end{table}The meaning of the labels in Table~\ref{table} is the following {\bf 1.} total execution time in seconds, {\bf 2.} number of integration 
	steps, {\bf 3.} if existence of a fixed point was proved, 
	{\bf 4.} if the fixed point is locally attracting, {\bf 5.} if the fixed point is attracting globally. Order of the Taylor method was $6$, time 
	step length was $0.005$ in all cases.
	
	For each case we fixed the radius of the energy absorbing ball $E_0$ and chose at random a forcing $f(x)$ which satisfies 
	$\frac{E(\{f_k\})}{\nu^2}=E_0$. 
	The forcing $f(x)$ was defined by a finite number of modes $\{f_k\}_{|k|\leq m}$. We added to each forcing mode $f_k$ the uniform 
	perturbation $[f_\varepsilon]:=[-\varepsilon,\varepsilon]\times[-\varepsilon,\varepsilon]$ (the parameter $\varepsilon$ is also provided
	in Table~\ref{table}) in order to perform simultaneously 
	a proof for a ball of functions. 
	\paragraph{}We would like to stress the fact that the provided cases are only examples and our program can attempt to prove any case.
	The package with the program along with the instruction and all the data from the proofs is available \cite{Package}.	
\section{Algorithm for constructing an absorbing set for large \gps}
	\label{sec:scb}
	\paragraph{}The goal of this section is to present an algorithm for constructing a set $V\oplus\mathnormal{\Theta}\subset\overline{H}$, forming 
	self-consistent bounds for \eqref{eq:burgers_infinite1} such that $V\oplus\mathnormal{\Theta}$ is an absorbing set for large \gps. 
	It is important to require that $V\oplus\mathnormal{\Theta}$ forms self-consistent bounds for \eqref{eq:burgers_infinite1} because in 
	Algorithm~\ref{alg:main} $V\oplus\mathnormal{\Theta}$ is integrated forward in time to verify that any solution in 
	$V\oplus\mathnormal{\Theta}$ after a finite time enters a trapping region.
	
	To support our claim that the constructed $V\oplus\mathnormal{\Theta}$ is in fact an absorbing set, in the following description we argue each estimate. 
	We drop the indication of Galerkin projections and times. For the precise meaning, the reader is referred to the proof of Lemma~\ref{lem:absorbingSet}. 
		\paragraph{Notation} $\sq{r}:=[-r,r]\times[-r,r]\subset\mathbb{R}^2,\quad \B{r}:=\overline{B(0,r)}\subset\mathbb{R}^2$.
		\paragraph{Input data} \begin{itemize} 
			\item $\nu>0,\quad M>m>0$ defining the dimensions of self-consistent bounds as in Definition~\ref{def:tail}, $\alpha\in\mathbb{R}$, 
			\item $\left\{[f_k]\right\}_{0<|k|\leq m}$ set of forcing modes perturbed by a uniform and constant perturbation $[f_\varepsilon]$, i.e. 
			$[f_k]=f_k+[f_\varepsilon]$ for $0<|k|\leq m$ and $[f_k]=0$ for $|k|>m$, $[f_0]=0$,
			\item $E_0$, where $E_0=\max_{\{f_k\}\in\left\{[f_k]\right\}}{\frac{E(\{f_k\})}{\nu^2}}$.\end{itemize}
		\paragraph{Output data} $V\oplus\mathnormal{\Theta}\subset\subspaceH$ forming self-consistent bounds for \eqref{eq:burgers_infinite1}.
		\paragraph{begin}		
		\paragraph{Initialization}$\hat{E}:=1.01\cdot \left(E_0+\alpha^2\right)$, $\hat{\varepsilon}:=10^{-15}$. 
		\paragraph{I Step} 
		\begin{itemize}
			\item For $0<|k|\leq M$ set 
			\begin{equation*}
				(V\oplus\mathnormal{\Theta})_k:=\sq{\frac{\hat{\varepsilon}+\left(\jednadruga\hat{E}+\max_{0<|k|\leq m}{\frac{\left|[f_k]\right|}{|k|}}\right)/\nu}{|k|}}.
			\end{equation*}
			\item For $|k|>M$ set 
			\begin{equation*}
				(V\oplus\mathnormal{\Theta})_k:=\B{\frac{\hat{\varepsilon}+\left(\jednadruga\hat{E}+\max_{0<|k|\leq m}{\frac{\left|[f_k]\right|}{|k|}}\right)/\nu}{|k|}}.
			\end{equation*}
		\end{itemize}
		\emph{Initial data is the absorbing ball of radius }$\hat{E}$\emph{, then by Lemma~\ref{lem:akbk2} combined with 
		Lemma~\ref{lem:energyEstimate} after a finite time the coefficients }$\{a_k\}$\emph{ satisfy}
		\begin{equation*}
			|a_k|\leq\frac{\hat{\varepsilon}+\left(\jednadruga\hat{E}+\max_{0<|k|\leq m}{\frac{|f_k|}{|k|}}\right)/\nu}{|k|}=:\frac{C}{|k|},\quad |k|>M.
		\end{equation*}
		\paragraph{II Step}		
		\begin{itemize}
		\item For $0<|k|\leq M$ calculate  
			\begin{equation*}
				b_{k,j}^-:=\frac{\left(-C\sqrt{\hat{E}}D+\frac{f^-_{k,j}}{|k|^\jednadruga}\right)/\nu}{|k|^\frac{3}{2}},\quad
				b_{k,j}^+:=\frac{\left(C\sqrt{\hat{E}}D+\frac{f^+_{k,j}}{|k|^\jednadruga}\right)/\nu}{|k|^\frac{3}{2}},\quad j=1,2.
			\end{equation*}
			\textit{Initial data is the set }$V\oplus\Theta$\emph{ from I Step, then the following estimate due to Lemma~\ref{lem:estimateNk} is
			used}
			\begin{equation*}
				|N_k(V\oplus\Theta)|\leq\frac{C\sqrt{\hat{E}}D}{|k|^{-\frac{1}{2}}},
			\end{equation*}
			\emph{where $C$ is defined in I Step.} 
		\item For $0<|k|\leq M$ set 
			\begin{equation*}
				(V\oplus\mathnormal{\Theta})_k:=
				[b_{k,1}^-,b_{k,1}^+]\times[b_{k,2}^-,b_{k,2}^+]+[-\hat{\varepsilon},\hat{\varepsilon}]^2.
			\end{equation*}	
			\emph{This is a refinement step. Using the data from I Step, a new value of }$V\oplus\Theta$\emph{ is defined.
			By Lemma~\ref{lem:akbk1} and Lemma~\ref{lem:refinementOfAbsorbingSet}
			after a finite time the coefficients }$\{a_k\}$\emph{ satisfy}
			\begin{equation*}
				a_k\in[b_{k,1}^-,b_{k,1}^+]\times[b_{k,2}^-,b_{k,2}^+]+[-\hat{\varepsilon},\hat{\varepsilon}]^2,
				\quad 0<|k|\leq M.
			\end{equation*}
		\item For $|k|>M$ set
			\begin{equation*}
				(V\oplus\mathnormal{\Theta})_k:=\B{\frac{\hat{\varepsilon}+\left(C\sqrt{\hat{E}}D\right)/\nu}{|k|^\frac{3}{2}}}.			
			\end{equation*}
			\emph{This is a refinement step. Using the data from I Step, a new value of }$V\oplus\Theta$\emph{ is defined.
			By Lemma~\ref{lem:akbk2} after a finite time the coefficients }$\{a_k\}$\emph{ satisfy}
			\begin{equation*}
				|a_k|\leq\frac{\hat{\varepsilon}+\left(C\sqrt{\hat{E}}D\right)/\nu}{|k|^\frac{3}{2}},\quad |k|>M.
			\end{equation*}
			\emph{Observe that $[f_k]=0$ for $|k|>M$ and Lemma~\ref{lem:akbk2} is used with $M_1=M$.}
		\end{itemize}
		\paragraph{III Step} Iterate the refinement, until $V\oplus\mathnormal{\Theta}$ forms self-consistent bounds for \eqref{eq:burgers_infinite1}, 
		as the stopping criterion use the condition $s(\mathnormal{\Theta})>d+p$, where $d$ and $p$ are from \eqref{eq:dPDEassumptions} and 
		$s(\mathnormal{\Theta})$ is the order of polynomial decay of the tail 
		$\mathnormal{\Theta}=\Pi_{|k|>M}{\overline{B\left(0,\frac{C(\mathnormal{\Theta})}{|k|^{s(\mathnormal{\Theta})}}\right)}}$. To calculate 
		$[b_{k,1}^-,b_{k,1}^+]\times[b_{k,2}^-,b_{k,2}^+]$ and $C(b)$ use the estimates derived in \cite{SuppMat}, in principle giving 
		much sharper bounds than energy-like estimates used in the previous steps. 
		
		\textit{Every such refinement generates 
		bounds that are reached by the solutions after a finite time. Moreover, to see that the procedure will stop, note that at each iteration 
		the \decay\ $s(\mathnormal{\Theta})$ is increased by $1$. Using the formulas derived in \cite{SuppMat} a bound such that $|N_k|\leq\frac{D}{|k|^{s(N)}}$ 
		is received, where $s(N)=s(\mathnormal{\Theta})-1$ and, therefore, $s_{new}(\mathnormal{\Theta})=s(b)=s(N)+2=s(\mathnormal{\Theta})+1$. 
		Finally, as soon as $s_{new}(\mathnormal{\Theta})>d+p$, stop.}
		\paragraph{end}
\section{Rigorous integration forward in time}
	By \textit{rigorous numerics} we mean algorithms for estimating solutions of differential equations that operate on sets and produce sets that 
	always contain an exact solution. Rigorous numerics for ODEs is a well established and analysed topic. There are a few algorithms
	that offer reliable computations of the solution trajectories for ODEs which are based on interval arithmetic. The approach used in this
	paper is based on the Lohner algorithm, presented in \cite{Lo} , see also \cite{ZLo}. It has made possible to prove many facts concerning
	the dynamics of certain ODEs, e.g. the Rossler equation, the Lorenz equation or the restricted n-body problem (see \cite{ZLo}, \cite{KZ} 
	and references therein). In the context of rigorous integration of ODEs we consider an abstract Cauchy problem
	\begin{equation}
		\left\{\begin{array}{ccc}
			\dot{x}(t)&=&f(x(t)),\\
			x(0)&=&x_0.
		\end{array}\right.
		\label{cauchy}
	\end{equation}
	$x\colon[0, \mathcal{T})\to\mathbb{R}^n$, $f\colon\mathbb{R}^n\to\mathbb{R}^n$, $f\in C^\infty$. 
	The goal of a rigorous ODEs solver is to find a set $\mathbf{x_k}\subset\mathbb{R}^n$ compact and connected such that
	\begin{equation}
		\varphi(t_k, \mathbf{x_0})\subset \mathbf{x_k},
	\end{equation} 
	$t_k\in[0, \mathcal{T}),\quad \mathbf{x_0}\subset\mathbb{R}^n$. By $\varphi(t_k, x_0)$ we denote the solution of \eqref{cauchy} at 
	the time $t_k$ with initial condition $x_0\in\mathbb{R}^n$, and therefore $\varphi(t_k, \mathbf{x_0})$ denotes the set of all the values which 
	are attained at the time $t_k$ by any solution of \eqref{cauchy} with the initial condition in $\mathbf{x_0}$.
	\paragraph{\textit{Notation}} We denote by $[x]$ an \textit{interval set} $[x]\subset\mathbb{R}^n$, 
	$[x]=\Pi_{k=1}^n[x^-_k, x^+_k]$, $[x^-_k, x^+_k]\subset\mathbb{R},\ -\infty<x_k^-\leq x_k^+<\infty$,
	$\midI{[x]}$ is the middle of an interval set $[x]$ and $\restI{[x]}$ is the rest, i.e. $[x]=\midI{[x]}+\restI{[x]}$.
	
	There are some subtle issues regarding intervals and set representation in the Lohner algorithm, which are discussed e.g. in 
	\cite{ZLo}. Let us only mention that it is highly ineffective to use the interval set representation explicitly $\Pi{[a^-_k, a^+_k]}$ because 
	it leads to the so-called \textit{wrapping effect} \cite{ZLo}, large over-estimates appear and prevents us from integrating over a longer 
	time interval. In order to avoid those problems we do not use interval sets explicitly, but to represent sets in a suitable coordinate system 
	we use the \textit{doubleton} representation of sets \cite{ZLo}
	\begin{equation}
		\label{eq:doubletonRepresentation}
		[x_k]+B_k\cdot [r_k]+C_k\cdot [r_0],
	\end{equation}
	where $B_k$ and $C_k$ are matrices representing a coordinate systems, $[x_k]$ is an interval set, likely a single point, $[r_k]$ is an interval set 
	that represents local errors that arise during integration, $[r_0]$ is an interval set that represents the error at the beginning (the diameter of a set 
	at the beginning).
	
	We stress the fact that we are interested in rigorous numerics for dPDEs, we develop main ideas in the following sections.
\subsection{Algorithm for integrating rigorously dPDEs}
	\label{sec:rigorousIntegration}
	In context of dPDEs we have to solve the following infinite system of ODEs
	\begin{equation}
		\label{eq:projection}
		\left\{
		\begin{array}{l}
			\frac{d x}{d t}=P_mF(x+y),\\
			\frac{d y}{d t}=Q_mF(x+y),\\
		\end{array}
		\right.
	\end{equation} 
	$x\in X_m$, $y\in Y_m$.
	
	Following \cite{KZ}, \cite{Z3} we will get estimates for \eqref{eq:projection} by considering the following differential 
	inclusion
	\begin{equation}
		\frac{d x}{d t}(t)\in P_mF(x(t))+\delta,
		\label{eq:delta}
	\end{equation}
	where $\delta\subset X_m$ describes influence of $y$ onto $P_mF(x+y)$. We call 
	\begin{equation}
		\label{eq:galerkinProjection}
		\frac{d x}{d t}=P_mF(x)
	\end{equation}
	the $m$ dimensional Galerkin projection of \eqref{eq:projection}, where $m>0$.
	
	We also consider a Cauchy problem, with $a\in X_m$, $x_0\in X_m$
	\begin{equation}
		\left\{\begin{array}{l}
			\frac{d x}{d t}(t)= P_mF(x(t))+a,\\
			x(0)=x_0.	
		\end{array}\right.
		\label{eq:associatedCauchyProblem}
	\end{equation}
	\paragraph{}Let $d_{X_m}$, $d_{Y_m}$ dimensions associated with $X_m$ and $Y_m$ respectively. From now on we switch to a more concrete setting, 
	which is 
	\begin{equation*}
		X_m:=\mathbb{R}^{d_{X_m}}\text{ and }Y_m:=\mathbb{R}^{d_{Y_m}},\ d_{X_m}<\infty,\ d_{Y_m}=\infty.
	\end{equation*}	
\providecommand{\spaceXm}{X_m}
\providecommand{\spaceYm}{Y_m}		
	In this section we assume that the solutions of problems \eqref{eq:projection}, \eqref{eq:galerkinProjection} and 
	\eqref{eq:associatedCauchyProblem} are defined and unique and later we will prove it.
	\paragraph{\textit{Notation}} $T, T(0), T(t_1), T([0,h])\subset\spaceYm$ 
	are tails satisfying \eqref{eq:farTail}, in the context of tails, for notational purposes, the symbol $T(\cdot)$ is not used to denote a function of 
	time, but an enclosure for a tail at the provided time. By 
	\begin{itemize}
		\item $\overline{\varphi}^{m}\left(t, x_0,a\right)$ we denote the solution of \eqref{eq:associatedCauchyProblem} at a time $t>0$ 
		with $a\in\spaceXm$ and an initial condition $x_0\in\spaceXm$,
		\item $\varphi^{X_m}\left(t,x_0,y_0\right)$ we denote the solution of \eqref{eq:projection} at a time $t>0$, projected onto $X_m$
		with an initial condition $x_0\in\spaceXm$ and $y_0\in\spaceYm$,
		\item $\varphi^{m}\left([0,h], x_0, T\right)$ denotes a collection of all possible values of the solution of the inclusion
		$\frac{d x}{d t}\in P_mF(x+T)$ on the time interval $[0,h]$ with $T\subset\spaceYm$ and an initial condition 
		$x_0\in\spaceXm$.
	\end{itemize}	 	
	Below we present all steps of the algorithm needed to rigorously integrate \eqref{eq:projection}. Whereas \cite{KZ}, \cite{Z3} algorithm is given in
	an abstract setting, here we provide a detailed description of an algorithm designed for dPDEs exclusively.
	\paragraph{} 
	In Algorithm~\ref{alg:main} we present steps needed to calculate rigorous bounds for the solutions of \eqref{eq:projection} at $t_1=h$. 
	The main idea is to get estimates for the solutions of each Galerkin projection of \eqref{eq:projection} simultaneously. For the correctness
	proof of Algorithm~\ref{alg:main} we refer the reader to \cite{KZ} or \cite{Z3}. Note that Algorithm~\ref{alg:main} is a subcase of a general 
	algorithm, with the set $[W_y]\subset\spaceXm$ chosen to be the Galerkin projection error.		
		\paragraph{Input}
		\begin{itemize}
			\item a time step $h$,
			\item $[f_\varepsilon]\subset\spaceXm$, a constant forcing perturbation,
			\item $[x_0]\subset\spaceXm$, an initial finite part,
			\item $T(0)\subset\spaceYm$, an initial tail,\\
				$[x_0]\oplus T(0)\subset H$ forms self-consistent bounds for \eqref{eq:evolutionFourierModes}.
		\end{itemize}
		\paragraph{Output}
		\begin{itemize}
			\item $[x_{t_1}]\subset \spaceXm$ such that $\varphi^{X_m}(t_1, [x_0], T(0))\subset[x_{t_1}]$, enclosure for the finite part of the solutions 
				at the time $t_1$.
			\item $T(t_1)\subset \spaceYm$, an enclosure for the tail at the time $t_1$,\\
				$[x_{t_1}]\oplus T(t_1)\subset H$ forms self-consistent bounds for \eqref{eq:evolutionFourierModes}.
		\end{itemize}
		\paragraph{begin}
		\begin{enumerate}
			\item find \label{step:first} $T\subset \spaceYm$ such that $T([0, h])\subset T$ and $[W_2]\subset \spaceXm$ such that 
				$\varphi^{m}([0, h], [x_0], T)\subset [W_2]$.
				Enclosure for the tail on the whole time interval $[0,h]$ and the enclosure for the collection of solutions of the differential 
				inclusion respectively. $[W_2]\oplus T$ forms self-consistent bounds for \eqref{eq:evolutionFourierModes},
			\item calculate the \textit{Galerkin projection error} $\spaceXm\supset[W_y]:=\{P_mF(x+T)-P_mF(x)\,|\,x\in[W_2]\}_I$,
			\item do the selection $[W_y]\ni y_c:=\midI{[W_y]}$, 
			\item \label{step:deterministic} apply the \textit{$C^0$ Lohner algorithm} to solve the system of autonomous ODEs \eqref{eq:associatedCauchyProblem} 
				with $a=y_c$. The result is a rigorous enclosure for the solution 
				$[\overline{x_{t_1}}]\subset\spaceXm\colon\overline{\varphi}^{m}(t_1, [x_0], y_c)\subset [\overline{x_{t_1}}]$.
				As a mid-step the enclosure $[W_1]$ such that $\overline{\varphi}^{m}\left([0,h], [x_0], y_c\right)\subset[W_1]$ is calculated and returned.
				Refer to \cite{ZLo} for the details,	
			\item calculate \textit{the perturbations vector} $\spaceXm\supset[\delta]:=\left[y_c-[W_y]+[f_\varepsilon]\right]_I$,	
			\item initialize the single valued vector $\spaceXm\ni C_i:=\sup{\left|[\delta_i]\right|}$,
			\item compute the ``Jacobian'' matrix $\mathbb{R}^{d_{X_m}\times d_{X_m}}\ni J\colon 
			J_{ij}\geq\left\{\begin{array}{ll}\sup{\frac{\partial F_i}{\partial x_j}([W_2], y_c)}&\text{if }i=j\\\left|\sup{\frac{\partial F_i}{\partial x_j}([W_2],y_c)}\right|&\text{if }i\neq j\end{array}\right.$,				
			\item perform component-wise estimates in order to calculate the set $[\Delta]\subset \spaceXm$,
				 $D:=\int_0^h{e^{J(t_1-s)}C}\,ds$,  $[\Delta_i]:=[-D_i, D_i]$ for $i=1,\dots,d_{X_m}$, 				
			\item obtain the final rigorous bound $[x_{t_1}]\subset \spaceXm$ for the solution of a differential inclusion by combining 
				results from the previous steps 
				$\varphi^{X_m}(t_1, [x_0], T(0))\subset[x_{t_1}]=[\overline{x_{t_1}}]+[\Delta],\quad [x_0]\subset \spaceXm,\ T(0)\subset \spaceYm$,
			\item perform rearrangements into the doubleton representation,
			\item \label{step:second} compute $T(t_1)\subset\spaceYm$ such that $\varphi^{Y_m}(t_1, [x_0], T(0))\subset T(t_1)$.
		\end{enumerate}
		\paragraph{end}
	\paragraph{}\RestyleAlgo{plain}\begin{algorithm}\label{alg:main}\caption{The main algorithm}\end{algorithm}\RestyleAlgo{boxed}
	\begin{rem}
	Basing on the framework of Algorithm~\ref{alg:main} we have developed an algorithm which apparently has been better in tests, the improvement concerns
	Step~\ref{step:first} and Step~\ref{step:second} of Algorithm~\ref{alg:main}. 
	As the details are very technical we do not present them here. The interested reader can find a detailed presentation in Appendix~\ref{sec:improvement},
	whereas in Appendix~\ref{sec:pseudoCode} we included the pseudo-codes. We omitted all the remaining steps of Algorithm~\ref{alg:main} that have already 
	been described in previous works. To realize some of the elements we used the \cite{CAPD} package.			
	\end{rem}
\section{Algorithm for proving Theorem~\ref{thm:main}}
	\label{sec:cap}
		\paragraph{\textit{Notation}} By a capital letter we denote \textit{a single valued matrix}, e.g. $A$, by $[A]$ we denote an \textit{interval matrix}. The inverse 
		matrix of $A$ is denoted by $A^{-1}$, we use the symbol $[A^{-1}]$ to denote an interval matrix such that $[A^{-1}]\ni A^{-1}$. $[M]_I$ denotes
		an interval hull of a matrix $M$, we also use this notation in the context of vectors. 
		\begin{figure}[htbp]
      \centering
      \includegraphics[width=\textwidth]{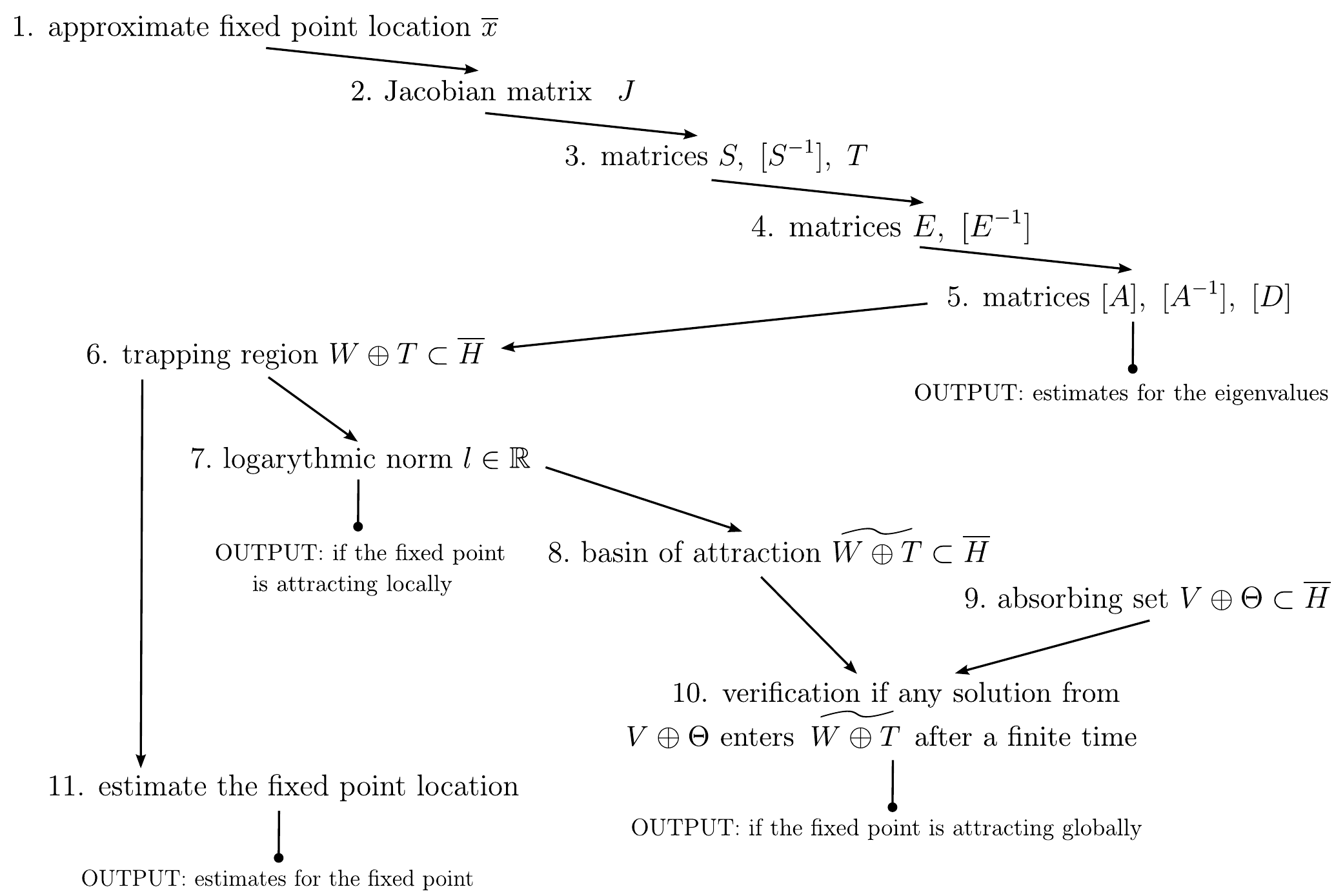}
      \caption{Flow diagram presenting steps of Algorithm for proving Theorem~\ref{thm:main}}
    \end{figure}		
		\paragraph{Input}
			\begin{itemize}
				\item $m>0$, an integer, the Galerkin projection \eqref{eq:symmetricGalerkinProjection} dimension,
				\item $[\nu_1,\nu_2]>0$, an interval of the viscosity constant values,
				\item $\alpha\in\mathbb{R}$, a constant value, equal to $\frac{1}{2\pi}{\int_{0}^{2\pi}{u_0(x)\,dx}}$,
				\item $s$, the \decay\ of coefficients that is required from the constructed bounds and trapping regions, have to be an integer 
					satisfying $s\geq 4$,
				\item order and the time step of the Taylor method used by the $C^0$ Lohner algorithm,
				\item set of $2\pi$ periodic forcing functions $f(x)$ for \eqref{eq:burgers}, defined by a finite number of modes $\{f_k\}_{0<|k|\leq m}$ 
					and a uniform and constant perturbation $[f_\varepsilon]=[-\varepsilon, \varepsilon]\times[-\varepsilon,\varepsilon]$.
			\end{itemize}
		\paragraph{Output}	
			\begin{itemize}
				\item $\overline{x}$, an approximate fixed point for \eqref{eq:burgers_infinite1},
				\item $J\approx dP_mF(\overline{x})$, an approximate Jacobian matrix at $\overline{x}$,
				\item $[A]$ and $[A^{-1}]$ interval matrices reducing $[dP_mF(\overline{x})]_I$ to an almost diagonal matrix $[D]$ - with dominating 
					blocks on the diagonal,
				\item $[D]=\left[\left[A\right]\cdot [dP_mF(\overline{x})]_I\cdot[A^{-1}]\right]_I$, almost diagonal form of the Jacobian matrix used 
					to estimate the eigenvalues of $dP_mF(\overline{x})$,
				\item $W\oplus T\subset\subspaceH$ and $\widetilde{W\oplus T}\subset\subspaceH$, trapping regions for \eqref{eq:burgers_infinite1} $W\oplus T\subset\widetilde{W\oplus T}$,
				\item $l$, upper bound of the logarithmic norm \eqref{eq:logarithmic} on $\widetilde{W\oplus T}$,
				\item $V\oplus \mathnormal{\Theta}\subset\subspaceH$, an absorbing set forming self-consistent bounds for \eqref{eq:burgers_infinite1},
				\item a rigorous bounds for the fixed point location,
				\item total time and integration steps needed to complete the proof.
			\end{itemize}
			\paragraph{begin}
			\begin{enumerate}
				\item find an approximate fixed point location $\overline{x}$ by non-rigorous integration of $\dot{x}=P_mF(x)$. Refine the provided 
					candidate $\overline{x}$ using \textit{the Newton method} iterations,
				\item calculate non-rigorously the Jacobian matrix, $J\approx dP_mF(\overline{x})$ (use $\nu_1$ as the viscosity constant in both 
					steps),
				\item calculate non-rigorously an approximate orthogonal matrix $S$ used for reducing $J$ to an approximate upper triangular matrix $T$
					(with 1x1 and 2x2 blocks on the diagonal). Use \textit{the QR algorithm} with multiple shifts to find $S$. Then find a rigorous inverse 
					$[S^{-1}]\colon S^{-1}\in[S^{-1}]$ using \emph{the Krawczyk operator} \cite{N},
				\item calculate \textit{the eigenvectors} of $T$ to form a block upper triangular matrix $E$ that is used to further reduce $T$ to 
				  an almost diagonal matrix, then calculate a rigorous inverse matrix $[E^{-1}]\colon E^{-1}\in[E^{-1}]$ using the Krawczyk operator 
				  again,
				\item calculate $[A]:=[E\cdot S]_I$, $[A^{-1}]:=[[S^{-1}]\cdot[E^{-1}]]_I$ and $[D]:=\left[\left[A\right]\cdot [dP_mF(\overline{x})]_I\cdot[A^{-1}]\right]_I$
					, where $[D]$ is in an almost diagonal form, having blocks on the diagonal and negligible intervals as non-diagonal elements, suitable
					form to estimate the eigenvalues,
				\item \label{step:trappingRegion} find $W\oplus T\subset\subspaceH$ a trapping region in block coordinates that encloses $\overline{x}$. This step requires  
					$[A]$ and $[A^{-1}]$, the change of coordinates matrices calculated in the previous step. A detailed description of an algorithm 
					performing this task is provided by \cite{ZAKS},
				\item calculate $l$ an upper bound for the logarithmic norm on the set $[[A^{-1}]\cdot W]_I\oplus T$, for the details refer to \cite{ZAKS}.									 
					In case $l<0$ by Theorem~\ref{thm:8} claim that there exists a locally attracting fixed point. Observe that in this case $W\oplus T$
					is the basin of attraction of the fixed point found. One may be tempted to use the ``analytical'' trapping region, calculated in Section~\ref{sec:analyticalTR} 
					for that purpose, but this is an unfeasible goal in general as an analytical trapping region may simply be too large to include 
					it into the calculation process,	
				\item enlarge $W\oplus T$ and return the largest calculated self-consistent bounds 
					$\widetilde{W\oplus T}\subset\subspaceH$ such that $\widetilde{W\oplus T}$ is a trapping region, $l<0$ and 
					$W\oplus T\subset\widetilde{W\oplus T}$. By Theorem~\ref{thm:8} claim that the basin of attraction of the fixed point found is
					$\widetilde{W\oplus T}$,
				\item using the procedure from Section~\ref{sec:scb} calculate the absorbing set $V\oplus \mathnormal{\Theta}\subset\overline{H}$,
				\item integrate $V\oplus \mathnormal{\Theta}$ rigorously forward in time until 
					$\varphi\left(t, \left[\left[A\right]\cdot V\right]_I\oplus\mathnormal{\Theta}\right)\subset\widetilde{W\oplus T}$ 
					. If this step finishes successfully conclude that $\widetilde{W\oplus T}$ after a finite time contains any solution of the problem 
					\eqref{eq:burgers_infinite} with sufficiently smooth initial data and claim the existence of a globally attracting fixed point,
				\item translate $[[A^{-1}]\cdot W]_I\oplus T$ into the doubleton representation \eqref{eq:doubletonRepresentation} 
					and integrate it forward in time in order to estimate the fixed point location with a relatively high accuracy.
			\end{enumerate}
			\paragraph{end}			
		\begin{rem}
			All the trapping regions constructed in the main algorithm presented above are expressed in block coordinates. Where the block decomposition
			of $H$ is given by $H=\oplus_{(i)}H_{(i)}$, where for $(i)>m$ blocks are given by $H_{(i)}=<e_i>$, and $(i)=i$ in this case. 
			Whereas for $(i)\leq m$  
			each block $H_{(i)}$ is a two-dimensional eigenspace of $J$. In case of two dimensional blocks $(i)=(i_1,i_2)\in\mathbb{Z}^2$,
			the expression $(i)<m$ means that $i_j<m$ for $j=1,2$. 
			
			Therefore given a trapping region $W\oplus T\subset H$ the finite part $W$ has the following form
			\begin{equation*}
				W=\prod_{(i)}\left\{\begin{array}{ll}
				\overline{B}\left(0, r_i\right)&\text{, for }(i)\in\mathcal{I},\\
				{[a_i^-,\ a_i^+]}&\text{, for }(i)\notin\mathcal{I},
				\end{array}\right.
			\end{equation*}			
		where $\mathcal{I}=\left\{(i)\colon H_{(i)}\text{ is two dimensional eigenspace of }J\right\}$. 
		\end{rem}		
		\begin{rem}
			\label{rem:statement}
			In all the proofs presented in Table~\ref{table} from Section~\ref{sec:exampleTheorem} we have
			\begin{equation*}
			\mathcal{I}=\left\{\begin{array}{ll}
			\emptyset&\text{ when $\int_0^{2\pi}{u_0(x)\,dx}=0$},\\
			\{(i)\colon(i)\leq m\}&\text{ when $\int_0^{2\pi}{u_0(x)\,dx}\neq 0$}.
			\end{array}\right.
			\end{equation*}
			We have not been able to prove this rigorously.
		\end{rem}	
\section{Conclusion}
A method of proving the existence of globally attracting fixed points for a class of dissipative PDEs has been presented. A detailed case study of the 
viscous Burgers equation with a constant in time forcing function has been provided. All the computer program sources used are available online 
\cite{Package}. There are several paths for the future development of the presented method we would like to suggest. An option is, for instance, 
to apply a technique for splitting of sets in order to see what the largest domain approachable by this technique is. One may also consider working 
on proving the statement given in Remark~\ref{rem:statement}. Another very interesting possibility is application of the presented method 
to higher dimensional PDEs, such as the Navier-Stokes equation, and we will address this topic in our forthcoming papers.   

\appendix	
\section{Data from the example proof}
	\label{sec:proofData}
	The parameters were as follows $\nu\in\paperExampleNu$ (the whole interval was inserted), $\widehat{m}=3$, $a_0=0.5$. To present the following data 
	all the numbers were truncated, for more precise data we refer the reader to the package with data from proofs available \cite{Package}. 
	
	The change of coordinates
	\begin{equation*}
		\mbox{\scriptsize
			$\midI{[A]}=
\left[\begin{array}{@{\,}c @{\,} c @{\,} c @{\,} c @{\,} c @{\,} c@{\,}}-0.998 & 0.0623 & -8.86\cdot 10^{-3} & -2.28\cdot 10^{-3} & -7.28\cdot 10^{-4} & -6.15\cdot 10^{-3}\\
0.0509 & 0.816 & 1.39\cdot 10^{-3} & -7.02\cdot 10^{-3} & -5.02\cdot 10^{-3} & 6.59\cdot 10^{-4}\\
6.35\cdot 10^{-3} & 0.0175 & 1.83\cdot 10^{-4} & 8.79\cdot 10^{-4} & -0.863 & 0.505\\
0.0175 & -6.12\cdot 10^{-3} & 6.8\cdot 10^{-4} & -1.86\cdot 10^{-4} & 0.505 & 0.863\\
0.0114 & 0.0153 & -0.288 & 0.957 & -6.9\cdot 10^{-4} & 4.41\cdot 10^{-5}\\
-0.0139 & 9.46\cdot 10^{-3} & 0.957 & 0.288 & -1.21\cdot 10^{-4} & 3.88\cdot 10^{-4}\\
\end{array}\right].
		$}
	\end{equation*}
	
	The Jacobian matrix in almost diagonal form
	\begin{equation*}
	\mbox{\scriptsize
		$\begin{split}
		&\midI{\left[[A]\cdot [dF(\overline{x})]_I\cdot[A^{-1}]\right]_I}\\
		=&
\left[\begin{array}{@{\,}c @{\,} c @{\,} c @{\,} c @{\,} c @{\,} c@{\,}}-2.06 & 0.402 & 0.0558 & -1.98\cdot 10^{-3} & 0.0123 & 0.0985\\
-0.598 & -2.04 & -7.17\cdot 10^{-3} & -0.0545 & -0.0985 & 0.0123\\
0.109 & -3.97\cdot 10^{-3} & -8.2 & 1 & 1.33\cdot 10^{-3} & 5.19\cdot 10^{-3}\\
-0.0143 & -0.112 & -1 & -8.2 & -5.19\cdot 10^{-3} & 1.33\cdot 10^{-3}\\
-0.0369 & 0.295 & -2\cdot 10^{-3} & 7.78\cdot 10^{-3} & -18.4 & 1.5\\
-0.295 & -0.0369 & -7.78\cdot 10^{-3} & -2\cdot 10^{-3} & -1.5 & -18.4\\
\end{array}\right].
		\end{split}$}
	\end{equation*}
	
Note that the matrix $\left[[A]\cdot [dF(\overline{x})]_I\cdot[A^{-1}]\right]_I$ does not have negligible elements beyond the diagonal blocks. This is
because we have performed the calculations for all the values $\nu\in\paperExampleNu$ simultaneously. If we perform the same calculations for one 
particular value of $\nu$ we would get a thin matrix with intervals of diameter $\thicksim10^{-15}$. 
	
The approximate eigenvalues $\spect(J)\approx$
	\begin{equation*}
	\mbox{\scriptsize
		$	\begin{split}
			\approx\left(-2.00088+0.489685i, -2.00088-0.489685i, -17.9982+1.50012i, -17.9982-1.50012i, -8.00096+0.999759i, -8.00096-0.999759i\right).
			\end{split}
		$}
	\end{equation*}
	The logarithmic norm upper bound \eqref{eq:logarithmic} $l = \paperExampleL$.\\
	The trapping region expressed in canonical coordinates $[[A^{-1}]\cdot \widetilde{W]_I\oplus T}=$
	\begin{equation*}
	\mbox{\scriptsize
	$
\begin{array}{|c|c|c|}\hline\mathbf{k} & \mathbf{\re{a_k}} & \mathbf{\im{a_k}}\\\hline\hline
1 & 2.59365\cdot 10^{-3}+[-0.144158,0.144158] & -6.66462\cdot 10^{-4}+[-0.171969,0.171969]\\
2 & 0.0984977+[-9.55661,9.55661]10^{-2} & -0.0123068+[-9.64073,9.64073]10^{-2}\\
3 & 4.57814\cdot 10^{-3}+[-5.7441,5.7441]10^{-2} & 0.0551328+[-5.7053,5.7053]10^{-2}\\
4 & -2.88994\cdot 10^{-4}+[-5.90827,5.90827]10^{-3} & -1.14901\cdot 10^{-3}+[-5.57297,5.57297]10^{-3}\\
5 & 1.01516\cdot 10^{-3}+[-1.99265,1.99265]10^{-3} & -2.33225\cdot 10^{-4}+[-2.46885,2.46885]10^{-3}\\
6 & 2.22928\cdot 10^{-5}+[-7.83055,7.83055]10^{-4} & 2.53646\cdot 10^{-4}+[-6.63366,6.63366]10^{-4}\\
7 & -1.08421\cdot 10^{-5}+[-1.60583,1.60583]10^{-4} & -2.06754\cdot 10^{-5}+[-1.48048,1.48048]10^{-4}\\
8 & 9.3526\cdot 10^{-6}+[-5.09454,5.09454]10^{-5} & -3.26429\cdot 10^{-6}+[-5.69385,5.69385]10^{-5}\\
\geq 9 & \multicolumn{2}{|c|}{|a_k|\leq0.970056/k^{4}}\\\hline\end{array}
	$}
	\end{equation*}
	The absorbing set $V\oplus \mathnormal{\Theta}=$
	\begin{equation*}
	\mbox{\scriptsize
	$
\begin{array}{|c|c|c|}\hline\mathbf{k} & \mathbf{\re{a_k}} & \mathbf{\im{a_k}}\\\hline\hline
1 & 4.96368\cdot 10^{-3}+[-0.142913,0.142913] & -2.33252\cdot 10^{-3}+[-0.144686,0.144686]\\
2 & 0.0971252+[-5.3667,5.3667]10^{-2} & -0.0125347+[-5.2554,5.2554]10^{-2}\\
3 & 4.25602\cdot 10^{-3}+[-2.65075,2.65075]10^{-2} & 0.0542654+[-2.71581,2.71581]10^{-2}\\
4 & -2.69437\cdot 10^{-4}+[-1.41625,1.41625]10^{-2} & -1.31697\cdot 10^{-3}+[-1.35799,1.35799]10^{-2}\\
5 & 1.06171\cdot 10^{-3}+[-7.40709,7.40709]10^{-3} & -2.2071\cdot 10^{-4}+[-7.931,7.931]10^{-3}\\
6 & -4.23977\cdot 10^{-6}+[-5.02386,5.02386]10^{-3} & 2.65357\cdot 10^{-4}+[-4.90303,4.90303]10^{-3}\\
7 & -2.23332\cdot 10^{-5}+[-3.434,3.434]10^{-3} & -3.58771\cdot 10^{-5}+[-3.42765,3.42765]10^{-3}\\
8 & 1.33681\cdot 10^{-5}+[-2.49967,2.49967]10^{-3} & -7.06016\cdot 10^{-6}+[-2.50257,2.50257]10^{-3}\\
\geq 9 & \multicolumn{2}{|c|}{|a_k|\leq147.297/k^{4}}\\\hline\end{array}
	$}
	\end{equation*}
	The absorbing set is apparently larger than the trapping region, it has been necessary to integrate it rigorously forward in time.
	The Taylor method used in the $C^0$ Lohner algorithm was of order $6$ with time step $0.005$.  Total execution time was $\paperExampleExecutionTime$ 
	seconds, total number of integration steps needed to verify that $\varphi\left(V\oplus \mathnormal{\Theta}\right)\subset \widetilde{W\oplus T}$ 
	(having in mind that the sets are expressed in different coordinates) was $\paperExampleNrOfSteps$, therefore $\widehat{t}=\paperExampleTotalTime$.
\section{Improvement of Algorithm~\ref{alg:main}}
\label{sec:improvement}
\providecommand{\CTO}{C_{T(0)}}
\providecommand{\sTO}{s_{T(0)}}
\providecommand{\MTO}{M_{T(0)}}
\providecommand{\CT}{C_T}
\providecommand{\sT}{s_T}
\providecommand{\MT}{M_T}
\providecommand{\CN}{C_N}
\providecommand{\sN}{s_N}
\providecommand{\MN}{M_N}
\providecommand{\CB}{C_b}
\providecommand{\sB}{s_b}
\providecommand{\MB}{M_b}
\providecommand{\CG}{C_g}
\providecommand{\sG}{s_g}
\providecommand{\MG}{M_g}
\providecommand{\CTH}{C_{T(h)}}
\providecommand{\sTH}{s_{T(h)}}
\providecommand{\MTH}{M_{T(h)}}
\subsubsection{Step~\ref{step:first} of Algorithm~\ref{alg:main}. The main loop.}
\label{sec:mainLoop}		
\begin{definition}
	\label{def:polynomialBound}
		Let $W\subset H$, $W$ convex. We call $W$ {\rm the polynomial bound} if there exist numbers $M>0$, $C>0$, $s\geq 0$
		such that 
			\begin{equation}
				\max_{x\in W_k}{||x||}\leq\frac{C}{|k|^s},\ |k|>M. 
			\end{equation}
		To denote the polynomial bound we use the quadruple $(W, M, C, s)$.
	\end{definition}
	Basically, during step~\ref{step:first} of Algorithm~\ref{alg:main} we have to calculate $T\subset Y_m$ a good enclosure for the tail  
	during the whole time interval $[0,h]$, i.e. $T$ has to satisfy $T([0,h])\subset T$. 
	Apparently, the bounds for $T([0,h])$ can be calculated explicitly, due to the following monotonicity of the bounds formula
	\begin{equation}
	  \label{eq:monotonicityOfBounds}
	  T([0,h])_k\subset T(0)_k\cup g_k,\quad k\in\mathbb{Z}\setminus\{0\},
	\end{equation}
	where $g_k$ is the linear approximation of the solution defined in Definition~\ref{def:nonlinearTerms}, see \cite[Lemma~6.1]{Z3}. $T(0)$ 
	in the formula \eqref{eq:monotonicityOfBounds} is known as it is the initial condition, and the polynomial bounds enclosing 
	$g$ can be calculated in a finite number of steps. We describe an appropriate procedure in the following part.
	\begin{definition}
		\label{def:nonlinearTerms}
		Let $W\oplus T\subset H$ forms a self-consistent bounds for \eqref{eq:evolutionFourierModes}, $m>0$ be the Galerkin projection 
		dimension, $N_k$, $f_k$ and $\lambda_k$ appear on the right-hand side of \eqref{eq:evolutionFourierModes}, $f_0=0$, $f_k=0$ for $|k|>m$. 
		For $k\in\mathbb{Z}\setminus\{0\}$ and $i=1,\dots,d_1$ we define 
		\begin{subequations}			
		\begin{gather}
			N_{k,i}^{\pm}\colon N_{k,i}^{-}\leq N_{k,i}(W\oplus T)\leq N_{k,i}^{+},\label{eq:boundsNk}\\
			b_{k,i}^{\pm}:=\frac{N_{k,i}^{\pm}+f_{k,i}}{-\lambda_k},\label{eq:boundsbk}\\
			g_{k,i}^{\pm}:=\left(T(0)_{k,i}^\pm-b_{k,i}^\pm\right)e^{\lambda_k h}+b_{k,i}^\pm\label{eq:boundsgk},\\
			N_k:=\Pi_{i=1}^{d_1}{[N_{k,i}^-,N_{k,i}^+]},\quad b_k:=\Pi_{i=1}^{d_1}{[b_{k,i}^-,b_{k,i}^+]},\quad g_k:=\Pi_{i=1}^{d_1}{[g_{k,i}^-,g_{k,i}^+]}\nonumber.
		\end{gather}
		\end{subequations}
	\end{definition}
		Now the question is how to verify the relations \eqref{eq:monotonicityOfBounds} in a finite number of steps. In general, it is impossible.
		Apparently, in the setting studied here, when sets are represented by polynomial bounds defined in Definition~\ref{def:polynomialBound} 
		the relations \eqref{eq:monotonicityOfBounds} can be verified in a finite number of steps. Observe that the self-consistent bounds
		introduced in Section~\ref{sec:generalMethod} are in particular polynomial bounds.
								
		We present procedures dealing with $\neartail$ and $\fartail$ in Algorithm~\ref{alg:validateNearTail} and 
		Algorithm~\ref{alg:validateFarTail}, to be found in Appendix~\ref{sec:pseudoCode}, separately for better clarification. For the exact 
		meaning of the symbols refer to Definition~\ref{def:tail}. The crucial part in Step~\ref{step:first} of Algorithm~\ref{alg:main} is to verify if 
		$\fartail([0,h])\subset\fartail$ in a finite number of steps, where $\fartail$ is a candidate for the far tail. 
		
		Now, let us present the procedure. Our goal is to enclose the interval sets $g_k$ by a uniform polynomial bound. Once we have a uniform
		polynomial bound, denoted by $g$, the verification of $\fartail([0,h])\subset\fartail$ is straightforward, because of the property 
		\eqref{eq:monotonicityOfBounds}. Firstly, given a polynomial bound
		\begin{equation}
			\label{eq:polynomialBoundT}		
			(W\oplus T, \MT, \CT, \sT)
		\end{equation}		
		a polynomial bound 
		\begin{equation*}
			(N, \MT, \CN, \sN)\text{ such that }\Pi_{k\in\mathbb{Z}}{N_k}\subset N
		\end{equation*}
		is found. This task requires performing some tedious estimates and we do not present them here. We derived the required estimates for 
		a class of dPDEs including the viscous Burgers equation in \cite{SuppMat}. Generally, a polynomial bound satisfying $\sN=\sT-r$ is found. 
		Then we immediately obtain a polynomial bound
		\begin{equation}
			\label{eq:polyBdBk}
			(b, \MT, \CB, \sB)\text{ such that }\Pi_{k\in\mathbb{Z}}{b_k}\subset b,
		\end{equation}
		with $\CB=\frac{\CN}{V(M)}$, $V(M)=\inf{\left\{\nu(|k|)\colon|k|>M\right\}}$ and $\sB=\sN+p$. Finally, a polynomial bound 
		\begin{equation}
			\label{eq:polyBdGk}
			(g, \MT, \CG, \sG)\text{ such that }\Pi_{k\in\mathbb{Z}}{g_k}\subset g
		\end{equation}
		is obtained using the formulas as follows
		\begin{lemma}If $|k|>M$ then
			\label{lem:formulasG}
			\begin{equation}
				\label{eq:Gformula1}
				|g_k|\leq\frac{\CTO e^{\lambda_{k}h}\cdot|k|^{\sB-\sTO}-\CB(e^{\lambda_{k}h}-1)}{|k|^{\sB}}					
			\end{equation}
			and
			\begin{equation}
				\label{eq:Gformula2}
				|g_k|\leq\frac{\CTO\cdot e^{\lambda_{k_{max}}h}(k_{max})^{\sB-\sTO}+\CB}{|k|^{\sB}}=\colon\frac{\CG}{|k|^{\sG}}
			\end{equation}
			where $k_{max}$ is the $k$ for which function $e^{\lambda_kh}\cdot k^{r}$ attains its maximum.
		\end{lemma}	
		\paragraph{	\textit{Proof}} Maximum of $f(k)=e^{\lambda_k h}k^r$, with ${\mbox{dom}\,{f}}=\{k\colon |k|>M\}$ is reached at $k_{max}$, therefore 
		\eqref{eq:Gformula1} is estimated by \eqref{eq:Gformula2} for any $|k|\geq M$.\qed
		\paragraph{}Note that 
		\begin{equation}
			\label{eq:powers}
			\sG>\sT
		\end{equation}
		because $\sN=\sT-r$, $\sG=\sB=\sT-r+p$ and $p>r$. 						
		\paragraph{The main loop}
		\paragraph{Input}
		$\left([x_0]\oplus T(0), \MTO, \CTO, \sTO\right)$ a polynomial bound, $[x_0]\oplus T(0)\subset H$ forms self-consistent bounds 
		for \eqref{eq:evolutionFourierModes}.
		\paragraph{Output}
		$\left([W_2]\oplus T, \MT, \CT, \sT\right)$ a polynomial bound such that $T([0,h])\subset T$, $\varphi^m\left([0,h],[x_0],T\right)\subset [W_2]$
			and $[W_2]\oplus T\subset H$ forms self-consistent bounds for \eqref{eq:evolutionFourierModes}.
		\paragraph{begin}
		\begin{enumerate}
			\item Initialize $T:=T(0)$.
			\item Update $\fartail$ using \fun{findS} function.
			\item \fun{while} \fun{not} \textit{validated}			
			\begin{itemize}		
				\item 
					$[W_2]:=$ \fun{enclosure}($[x_0]$, $T$),
					calculate a rough-enclosure $[W_2]$ for the differential inclusion \eqref{eq:delta} using a current candidate for a
					tail enclosure $T$, after this step $\varphi^m\left([0,h],[x_0],T\right)\subset [W_2]$ holds,
				\item 
					calculate the polynomial bounds $(b, \MT, \CB, \sB)$ and $(g, \MT, \CG, \sG)$,					
				\item 
					\textit{validated} $:=$ \fun{validateTail}($T(0)$, $T$, $b$, $g$, $[W_2]$) (if $T$ was changed during this step \textit{validated}=\fun{false}).
			\end{itemize}
			\fun{end while}
		\end{enumerate}
		\paragraph{end}
		\begin{rem} In our algorithm the number $\MT$ in \eqref{eq:polynomialBoundT} is chosen adaptively in \fun{validateFarTail} and changes 
		from step to step.
		\end{rem}
		\paragraph{}\fun{enclosure} is \textit{the rough enclosure algorithm based on isolation}, designed for dPDEs presented in \cite{Z3}.
		
		We present a correctness proof of the \fun{validateNearTail} and \fun{validateFarTail} functions in the comments within the code listings from 
		Appendix~\ref{sec:pseudoCode}. By a correctness proof we show that a polynomial bound $T$, such that the condition  
		$T([0,h])_k\subset T_k$ holds for all $k\in\mathbb{Z}$, is returned by the algorithm whenever the algorithm stops.
		
		Now, we shall focus on explaining the main idea behind \fun{validateFarTail} and explain why we consider it an improvement of the existing 
		algorithm. Basically, when a $-\lambda_k$ in \eqref{eq:evolutionFourierModes} is small, the nonlinear part $N_k$ dominates the linear term. 
		However, there exists an index $\tilde{k}\in\mathbb{N}$ such that $-\lambda_k$ for $|k|>\tilde{k}$ becomes large enough to make the linear 
		part overtake the nonlinear part. The position of the threshold $\tilde{k}$ depends on the maximum order of the ``Laplacian'' that appears in 
		the linear part $L$ of \eqref{eq:evolutionFourierModes}, as well as on the order of the polynomial that appears in the nonlinear part. We 
		remark that the solution of the $m$-th Galerkin projection of \eqref{eq:evolutionFourierModes} with $m<\tilde{k}$ greatly differs from the 
		solution of the whole system \eqref{eq:evolutionFourierModes}.

		The aforementioned effects show that a proper choice of the Galerkin projection dimension $m$ (in our algorithm taken only once at the beginning) 
		and the number $\MT$ of the polynomial bound \eqref{eq:polynomialBoundT} (in our algorithm taken at each time step) is of critical importance and 
		has to be performed carefully. The application of a too small value may result in blow-ups and may prevent the completion of the calculations. 
		In the original algorithm from \cite{Z3} the number $\MT$ was fixed in advance.
		Then heuristic formulas were derived for the KS equation in order to predict if the tail validation function would finish successfully
		for a given $\MT$, $\sT$ and to guess the initial values of $\CT$ and $\sT$ in \eqref{eq:polynomialBoundT}, see \cite[Section~8]{Z3}.  
		We found the original approach insufficient for the purpose of 
		rigorously integrating PDEs that are the subject of our research (for example the Burgers or the Navier-Stokes equations). 
		A similar approach for the mentioned dPDEs is problematic and, especially in the case of lower viscosities, heuristic formulas cause 
	  performance issues and sometimes offer infeasible values, mainly due to the lower order of the ``Laplacian'' in the linear part.	
\subsubsection{Step~\ref{step:second} of Algorithm~\ref{alg:main}}
	\paragraph{Input}
		$\left([W_2]\oplus T, \MT, \CT, \sT\right)$, a polynomial bound such that $T([0,h])\subset T$,  
		$\varphi^m\left([0,h],[x_0],T\right)\subset [W_2]$ and $[W_2]\oplus T\subset H$ forms self-consistent bounds for \eqref{eq:evolutionFourierModes}.
	\paragraph{Output}
		$\left(T(h), \MTH, \CTH, \sTH\right)$, a polynomial bound.
	\paragraph{begin}
	\begin{enumerate}
		\item $\MTH:=\MT$, $T(h)$ inherits $M$ from the enclosure $T$,
		\item calculate the polynomial bound $(g, \MT, \CG, \sG)$,
		\item $T(h):=g$, $C_{T(h)}:=C_g$, $s_{T(h)}:=s_g$.
	\end{enumerate}
	\paragraph{end}	

\section{Validate tail function in pseudo-code}
\label{sec:pseudoCode}
	Here we present a pseudo-code of the functions \fun{validateNearTail} and \fun{validateFarTail} used in Section~\ref{sec:mainLoop}. First, we 
	present the internal representation of sets that was used in actual program, written in C++ programming language and available at \cite{Package}.
\paragraph{\textit{Data representation}}	
	\begin{itemize}
		\item \textit{double} is a floating point number of double precision in \textit{C++ programming language},
		\item \textit{interval} is $[a^-, a^+]\subset\mathbb{R}$ where $a^-$, $a^+$ are \textit{double} numbers. All arithmetic operations 
		on such intervals are rigorous and are performed using implementation of the CAPD library \cite{CAPD}. It is verified that the interval 
		arithmetic provides proper in mathematical sense results \cite{N},
		\item \textit{Vector} represents an interval set, a vector composed of \textit{intervals},
		\item \textit{PolyBd} is a structure used for representing a polynomial bound $(W, M, C, s)$. A given \textit{PolyBd} $V$ contains a 
		\textit{Vector} representing the finite part of $W\subset H$, an \textit{integer} representing $M$ denoted by $M(V)$ and two \textit{intervals} 
		representing $C$ and $s$ denoted by $C(V)$ and $s(V)$ respectively.
	\end{itemize}
	Below, in Algorithm~\ref{alg:validateNearTail} and Algorithm~\ref{alg:validateFarTail}, we present functions \fun{validateNearTail} and 
	\fun{validateFarTail} respectively along with correlated functions. Wherever $previous$ keyword appear the value from the previous step is used.\\
	\begin{function}[H]
		\TitleOfAlgo{predictM}
		\SetAlgoNoLine
		\KwIn{\textit{PolyBd} $T$, \textit{PolyBd} $g$}
		$L:=\left(\frac{C(T)}{C(g)}\right)^{s(T)-s(g)}$\;
		\Return $L$\;
	\end{function}	
	\begin{function}[H]
		\TitleOfAlgo{correctM}
		\SetAlgoNoLine
		\KwIn{\textit{PolyBd} $T$, \textit{PolyBd} $T(0)$, \textit{double} $L$}
		\tcp{function corrects current dimension $M$ of tails in two cases: value of $L$ is increasing}
		\If{$L>$previous $L$}{
			\lIf{$L$ is sufficiently small}{$M(T):=M(T(0)):=L_d$\;}
			\lElse{$M(T):=M(T(0)):=M$\;}
		}
		\tcc{and test if $M$ can be decreased, by checking if $L$ have established,
		by comparing approximation of current and previous $L$ up to the order $10^2$}
		\If{\fun{truncate}($L$, 2) $=$ \fun{truncate}(previous $L$, 2)}{
			$M(T):=M(T(0)):=L$\;
		}	
	\end{function}	
	\begin{function}[H]
	\TitleOfAlgo{findS}	
	\SetAlgoNoLine
	\KwIn{\textit{PolyBd} $T(0)$, \textit{PolyBd} $T$, \textit{Vector} $W_2$}
	\tcc{heuristic function, tries to find optimal $s(T)$ at each iteration of the main loop. By optimal $s(T)$ we mean largest possible value such 
	that a predicted $M$ is within desired range. We recall that we start with s(T)=s(T(0))}
	$PolyBd$ $g:=g(T(0), T, W_2)$\;
	$currentM:=$ \FuncSty{predictM}($T$, $g$)\;
	$potentialM:=$ \FuncSty{predictM}($T$, $g$ with decreased $s$)\;
	\While{$currentM$ out of desired range \KwSty{and} $s(T)>p+d$ \KwSty{and} $potentialM>2m$}{
		$currentM:=$ \FuncSty{predictM}($T$, $g$)\;
		$potentialM:=$ \FuncSty{predictM}($T$, $g$ with decreased $s$)\;
		$C(T):=C$(previous $T$)$\cdot\left(M+1\right)^{s(T)-s(\text{previous} T)}$\;
		$s(T):=s(T)-1$\;
	}
	\lIf{ $s(T)$ != previous $s(T)$}{
		
		\FuncSty{correctM}($currentM$)\;
	}
	\end{function}
	\begin{function}
	\TitleOfAlgo{update}
	\SetAlgoNoLine
	\KwIn{\textit{PolyBd} $T$, \textit{PolyBd} $T'$, \textit{set} $\mathcal{I}$}
	\For{$i\colon i\in\mathcal{I}$}{
		\If{$T'_i\nsubseteq T_i$}{
			\fun{calculate} \fun{new} $T_i\colon T'_i\subset$ \fun{new} $T_i$\;
			$T_i:=$\fun{new} $T_i$\;
		}
	}	
	\end{function}		
	\begin{algorithm}[H]
	\SetAlgoNoLine						
	\KwIn{\textit{PolyBd} $T(0)$, \textit{PolyBd} $T$, \textit{PolyBd} $g$, \textit{Vector} $W_2$}
	\KwOut{$bool$}
		\tcc{individually verify condition $T(0)_i\cup g_i\subset T_i$}
		$vector<bool> inflatesRe$\;
		$vector<bool> inflatesIm$\;
		\For{$k:=m+1,\dots,M$}{
		\If{!($\re{b_k}^+\leq\re{T(0)_k}^+$) \KwSty{and} !($\re{T_k}^+>\re{g_k}^+$)}{
			$\re{T_k^+}:=\re{g_k^+}$\;
			$inflateRe:=true$\;
		}
		\If{!($\re{b_k}^-\geq\re{T(0)_k}^-$) \KwSty{and} !($\re{T_k}^-<\re{g_k}^-$)}{
			$\re{T_k^-}:=\re{g_k^-}$\;
			$inflateRe:=true$\;
		}			
		\If{!($\im{b_k}^+\leq\im{T(0)_k}^+$) \KwSty{and} !($\im{T_k}^+>\im{g_k}^+$)}{
			$\im{T_k^+}:=\im{g_k^+}$\;
			$inflateIm:=true$\;
		}
		\If{!($\im{b_k}^-\geq\im{T(0)_k}^-$) \KwSty{and} !($\im{T_k}^-<\im{g_k}^-$)}{
			$\im{T_k}^-:=\im{g_k}^-$\;
			$inflateIm:=true$\;
		}
		\If{$inflateRe$}{				
			\fun{inflate}($\re{T_k}$, $1+c_{inflate}$)\;
			\For{ $j:=-c_{radius},\dots,c_{radius}$}{
				$inflatesRe[k+j]:=inflatesRe[k+j]+1+c_{inflate}/|j|$\;
			}
		}
		\If{$inflateIm$}{
			\fun{inflate}($\im{T_k}$, $1+c_{inflate}$)\;
			\For{ $j:=-c_{radius},\dots,c_{radius}$}{
				$inflatesIm[k+j]:=inflatesIm[k+j]+1+c_{inflate}/|j|$\;
			}
		}
		}
		\For{$k:=m+1,\dots,M$}{
			\If{$inflatesRe[k]>0$}{
				\fun{inflate}($\re{T_k}$, $inflatesRe[k]$)\;
			}
			\If{$inflatesIm[k]>0$}{
				\fun{inflate}($\im{T_k}$, $inflatesIm[k]$)\;
			}
		}
		\caption{validateNearTail function}
		\label{alg:validateNearTail}
	\end{algorithm}	
	\begin{algorithm}[H]
	{\scriptsize
	\caption{validateFarTail function}
	\label{alg:validateFarTail}
	\SetAlgoNoLine		
	\KwIn{\textit{PolyBd} $T(0)$, \textit{PolyBd} $T$, \textit{PolyBd} $b$, \textit{PolyBd} $g$, \textit{Vector} $W_2$}
	\KwOut{$bool$} 	 	
		$L:=\left(\frac{C(T)}{C(g)}\right)^{\frac{1}{s(T)-s(g)}}$; $L_2:=\left\lceil\left(\frac{C(b)}{C(T(0))}\right)^{\frac{1}{s(b)-s(T0)}} \right\rceil$\;		
		\paragraph{Case 1} $s(b)>s(T(0))$\\	
			\eIf(\tcp*[f]{in particular $T(0)_{M+1}\subset g_{M+1}\subset b_{M+1}$}){$T(0)_{M+1}\subset b_{M+1}$}{
				\lIf{$L_2<M+1$}{\fun{throw}($exception$)}\;
				\lIf{$T_{M+1}\varsubsetneq g_{M+1}$}{
					\fun{update}($T,g,\{M+1,M+2,\dots\}$)\;
				}
				\If{$L_2<\infty$}{
					\lIf{$T_{L_2}\varsubsetneq T(0)_{L_2}$}{\fun{update}($T,T(0),\{L_2,L_2+1,\dots\}$)\;}						
				}
				\tcc{If $L_2=\infty$ it is enough to validate $T_{M+1}$ only, because $s_g>s_T$, see \eqref{eq:powers}. If $L_2<\infty$, 
				$T_{M+1}$ is validated to cover the finite number of indices $\{M+1,\dots,L_2\}$
				and then validating $T_{L_2}$ covers the infinite rest $\{L_2,L_2+1,\dots\}$ due to $s(T)\leq s(T(0))$, see \fun{findS} function.}
				\lIf{$T$ was updated}{\fun{correctM}($T$, $T(0)$, $L$)\;}
			}(\tcp*[f]{$b_{M+1}\varsubsetneq T(0)_{M+1}$, in particular $b_{M+1}\varsubsetneq g_{M+1}\varsubsetneq T(0)_{M+1}$}){
				\lIf{$T_{M+1}\varsubsetneq T(0)_{M+1}$}{\fun{update}($T, T(0), \{M+1,M+2,\dots\}$)\;}
				\tcc{It is enough to validate $T_{M+1}$ only, because $s(T)\leq s(T(0))$ and $b_i\varsubsetneq T(0)_i$ for all $i>M$.}
			}
		\paragraph{Case 2} $s(b)=s(T(0))$\\
		\eIf{$b_{M+1}\subseteq T(0)_{M+1}$}{
			\lIf{$T_{M+1}\varsubsetneq T(0)_{M+1}$}{\fun{update}($T, T(0), \{M+1,M+2,\dots\}$)\;}
		}(\tcp*[f]{$T(0)_{M+1}\varsubsetneq b_{M+1}$}){
			\If{$T_{M+1}\varsubsetneq g_{M+1}$}{
				\fun{update}($T,g,\{M+1,M+2,\dots\}$)\;
				\fun{correctM}($T$, $T(0)$, $L$)\;
			}
		}
		\tcc{In both cases it is enough to validate $T_{M+1}$ because either $b_i\subset T(0)_i$ for all $i>M$ or $T(0)_i\subset b_i$ for all $i>M$ 
		and $s(T(0))=s(b)=s(g)\geq s(T)$, see \eqref{eq:powers}.}
		\paragraph{Case 3} $s(b)<s(T(0))$\\
		\eIf{$b_{M+1}\subset T(0)_{M+1}$}{
			\lIf{$L_2<M+1$}{\fun{throw}($exception$)\;}
			\lIf{$T_{M+1}\varsubsetneq T(0)_{M+1}$}{
				\fun{update}($T,T(0),\{M+1,M+2,\dots\}$)\;
			}
			\If{$L_2<\infty$}{
				\lIf{$T_{L_2}\varsubsetneq g_{L_2}$}{\fun{update}($T,g,\{L_2,L_2+1,\dots\}$)\;}
			}
			\tcc{If $L_2=\infty$ it is enough to validate $T_{M+1}$ only, because $s(T)\leq s(T(0))$. If $L_2<\infty$, 
			$T_{M+1}$ is validated to cover the finite number of indices $\{M+1,\dots,L_2\}$ and validating $T_{L_2}$ 
			covers the infinite rest $\{L_2,L_2+1,\dots\}$ due to $s_g>s_T$, see \eqref{eq:powers}.}
			\lIf{$T$ was updated}{\fun{correctM}($T$, $T(0)$, $L$)\;}
		}(\tcp*[f]{$T(0)_{M+1}\varsubsetneq b_{M+1}$}){
			\If{$T_{M+1}\varsubsetneq g_{M+1}$}{
				\fun{update}($T,g,\{M+1,M+2,\dots\}$)\;
				\fun{correctM}($T$, $T(0)$, $L$)\;
			}
			\tcc{It is enough to validate $T_{M+1}$ only, because $T(0)_i\varsubsetneq b_i$ for all $i>M$ and $s_g>s_T$, see \eqref{eq:powers}.}
		}
		\lIf{$T$ was updated}{\fun{return} false\;}\lElse{\fun{return} true\;}
	}	
	\end{algorithm}	
\end{document}